\documentclass[a4paper, 11pt, twoside]{article}

\usepackage{appendix}
\usepackage{amsmath,amsthm,empheq}
\usepackage{amssymb,latexsym}
\usepackage{mathrsfs}
\usepackage{enumerate}

\usepackage[colorlinks,citecolor=blue]{hyperref}
\usepackage[numbers,sort&compress]{natbib}
\usepackage{indentfirst}

\DeclareMathAlphabet{\mathpzc}{OT1}{pzc}{m}{it}

\newtheorem{theorem}{Theorem}[section]

\newtheorem{lemma}{Lemma}[section]
\newtheorem{proposition}{Proposition}[section]

\newtheorem{remark}{Remark}[section]

\newtheorem{example}{Example}[section]

\theoremstyle{definition} \theoremstyle{remark}
\numberwithin{equation}{section}
\allowdisplaybreaks

\begin{document}
	\markboth{X. Huang, L. Peng, J.C. Pozo, Y. Zhou}{Fully nonlocal telegraph equations}
	
	\date{}
	\baselineskip 0.22in
	\title{{\bf Random data Cauchy theory for fully nonlocal telegraph equations }}
	
	\author{Xi Huang$^{1}$, Li Peng$^{1}$\thanks{\footnotesize {Corresponding author.} }, Juan Carlos Pozo$^{2}$, Yong Zhou$^{1,3}$ \\[1.8mm]
		\footnotesize  {$^1$Faculty of Mathematics and Computational Science, Xiangtan University,}\\
		\footnotesize  {Hunan 411105, China}\\
		\footnotesize  {$^2$Acad\'{e}mico Instituto de Ciencias de la Ingenier\'{\i}a,
Universidad de O'Higgins, Rancagua, Chile}\\
\footnotesize  {$^3$Faculty of Information Technology, Macau University of Science and Technology,}\\
		\footnotesize  {Macau 999078, China}
	}
	
	\maketitle
	
	\begin{abstract}
		We consider the random Cauchy problem for the fully nonlocal telegraph equation of power type with the general $(\mathcal{PC}^{\ast})$ type kernel $(a,b)$. This equation can effectively characterize high-frequency signal transmission in small-scale systems. We establish a new completely positive kernel induced by $b$ (see Appendix \refeq{app b}) and derive two novel solution operators by using the relaxation functions associated with the new kernel,
which are closely related to the operators $\cos(\theta(-\Delta)^{\frac{\beta}{4}} )$ and $(-\Delta)^{-\frac{\beta}{4} }\sin(\theta(-\Delta)^{\frac{\beta}{4}} )$ for $\beta\in(1,2]$. These operators enable, for the first time, the derivation of mixed-norm $L_t^qL_x^{p'}$ estimates for the novel solution operators. Next, utilizing probabilistic randomization methods, we establish the average effects, the local existence and uniqueness for a large set of initial data $u^\omega \in L^{2}(\Omega, H^{s,p}(\mathbb R^3))$ ($p\in (1,2)$) while also obtaining probabilistic estimates for local existence under randomized initial conditions. The results reveal a critical phenomenon in the temporal regularity of the solution regarding the regularity index $s$ of the initial data $u^\omega$. \\[2mm]
	\end{abstract}
	
	\tableofcontents
	
	\section{Introduction}
It is well known that the classical telegraph equations adequately describe electromagnetic wave propagation along transmission lines for low-frequency signal transmission in large-scale systems with good conductors. However, it becomes necessary to re-examine phenomena such as charge accumulation along lines and memory effects in polarization and magnetization processes when addressing high-frequency signal transmission in small-scale systems. To model non-Maxwellian wave propagation in complex media, Pozo and Vergara \cite{J. C. Pozo 20} began with the following transmission line balance equations incorporating constitutive effects
	\begin{align*}
		&\partial_{x} V(t,x)+\frac{D h}{W} \partial_{t}(\mathrm{d} \varepsilon_1\ast I(\cdot,x))(t)+\frac{h}{A}\left(\mathrm{d} \overline{\varepsilon}_{R} \ast I(\cdot,x)\right)(t)=0,\\
		&\partial_{x} I(t, x)+\frac{A}{h} \partial_{t}(\mathrm{d} \varepsilon_2 \ast V(\cdot,x))(t)+\frac{A}{h}\left(\mathrm{d} \varepsilon_{G} \ast V(\cdot,x)\right)(t)=0,
	\end{align*}
	where $I$ and $V$ denote total tension and total current, respectively, $R$ resistance, $G$ the leak-conductance, $A$ cross-sectional area, $D$ longitudinal length, $W$ lateral width and $h$ medium thickness. $\varepsilon_1, \varepsilon_2$, $\varepsilon_G, \overline{\varepsilon}_R \in BV_{loc}(\mathbb R_+)$ are given material functions. For a function $\varepsilon \in BV_{loc}(\mathbb R_+)$,  $(\mathrm{d}\varepsilon\ast w)(t) = \int_0^t w(t-\zeta)\, \mathrm{d}\varepsilon(\zeta)$. Then the above relation can be derived as the generalized time nonlocal telegraph equations
	\begin{align*}
		\mu_1\partial_{t}^{2}(\mathrm{d} \varepsilon_1 \ast \mathrm{d} \varepsilon_2 \ast w)+\partial_{t}\left(\left[\mathrm{d} \overline{\varepsilon}_{R} \ast \mathrm{d} \varepsilon_2+\mu_1 d \varepsilon_1 \ast \mathrm{d} \varepsilon_{G}\right] \ast w\right)+ \mathrm{d} \overline{\varepsilon}_{R} \ast \mathrm{d} \varepsilon_{G} \ast w-\partial_{x}^{2} w=0,
	\end{align*}
	where $\mu_1 := \frac{AD}{W}$ and $w$ could be defined as $I$ or $V$. A simplified and commonly studied model that captures key features of the above system is given by
	\begin{align}\label{1.1}
		\partial_{t}^{2}(a_1 \ast a_1 \ast w(\cdot,x))(t) +\gamma\partial_{t}\left(a_1 \ast w(\cdot,x)\right)(t)-\partial_{x}^{2} w=0, \quad x \in \mathbb R,
	\end{align}
	where $\gamma>0$ and $a_1 \in L_{loc}^1(\mathbb R_+)$ satisfy the following condition
	\begin{enumerate}
		\item[$(\mathcal{PC})$] $a_1 \in L^1_{loc}(\mathbb R_+)$ is nonnegative and nonincreasing, and there exists $b \in L_{\text{loc}}^1(\mathbb{R}_+)$ such that $(a_1 \ast b)(t) = 1$ for $t>0$.
	\end{enumerate}
	In particular, when taking $a_1=g_{1-\alpha}$ with $\alpha \in (0,1)$, \eqref{1.1} reduces to the time fractional telegraph equation, which  has attracted significant attention due to its ability to nonlocal phenomena in transmission line media, as demonstrated in \cite{E. Orsingher, R. Figueiredo Camargo, R.C. Cascaval, R. Ashurov}.
	
	The study of $(\mathcal{PC})$ class kernels originally appeared in the context of nonlocal-in-time subdiffusion problems. This research direction proved particularly fruitful because such equations admit direct reformulations as Volterra integral equations with completely positive kernels, leading to several breakthrough results in the field. Vergara and Zacher \cite{V. Vergara} gave sharp decay estimates of solution for the equation
	\begin{align*}
		\partial_t a_1\ast(w-w_0) -\mathrm{div}(A(t,x)Dw) = 0
	\end{align*}
	in bounded domains with homogeneous Dirichlet boundary conditions, where the coefficient matrix $A$ satisfies measurability, boundedness, and the uniform parabolicity condition. Their proofs relied on energy methods and a new inequality for integro-differential operators. Kemppainen \cite{J. Kemppainen} et al. investigated the time-decay estimates for solution to
	\begin{align*}
		\begin{cases}
			\partial_t a_1\ast(w-w_0) -\Delta w = 0, \quad t>0, \ x \in \mathbb R^N,\\
			w_{t=0} = w_0, \quad x \in \mathbb R^N.
		\end{cases}
	\end{align*}
	Based on Fourier multiplier methods and properties of relaxation functions, they obtained optimal $L^p$-decay rates, and then they claimed that the decay profile on solutions presents a critical dimension
	phenomenon. More research on nonlocal subdiffusion problems can be found in \cite{J. C. Pozo 24, J. Kemppainen 17, J.C. Pozo}.
	
	Current research focuses on investigating the connection between fundamental solutions of time nonlocal telegraph equations with convolution kernel $a_1 \in (\mathcal{PC})$ and stochastic processes. Pozo and Vergara \cite{J. C. Pozo 20} established that the fundamental solution of \eqref{1.1} corresponds to the probability density function of a specific stochastic process $X(t)$. Furthermore, employing the Karamata-Feller Tauberian theory (\cite[Chapter XIII]{W. Feller}), they characterized the asymptotic behavior of the process variance across different temporal scales. Subsequently, Alegr\'{i}a and  Pozo \cite{F. Alegria} introduced a systematic technique to generate examples where the variance $\text{Var}X(t)$ grows sublinearly or logarithmically in time, contrasting with classical diffusive scaling. Furthermore,  Thang \cite{N. N. Thang} established that for any given completely monotone ultraslow kernel $a_1$ in \eqref{1.1}, there exists an induced kernel such that the mean squared displacement of the associated stochastic process for the time nonlocal telegraph equation exhibits logarithmic growth. Recently, Alegr\'{i}a and Pozo \cite{F. Alegria 23} considered a non-Markovian telegraph processes $X_k(t)$
	whose probability density function solves the following problem
	\begin{align*}
		\begin{cases}
			\partial_{t}^{2}(a_1 \ast a_1 \ast (w(\cdot,x)-\delta_0))(t) +\gamma\partial_{t}\left(a_1 \ast (w(\cdot,x)-\delta_0) \right)(t)-\Delta w=0, \quad t>0, \ x \in \mathbb R^N,\\
			w(0,x) = \delta_{0}(x), \ \partial_t w(0,x) = 0, \quad x \in \mathbb R^N,
		\end{cases}
	\end{align*}
	where $\delta_0(x)$ is the Dirac distribution.
	They rigorously established that the moments of $X_k(t)$ satisfies the Carleman determinacy condition, while further demonstrating that its probability distribution admits a representation through subordination of the classical telegraph process $T(t)$ by the random time change $|W(t)|$ related to a time nonlocal wave equation.
	
	In summary, research progress on the generalized nonlocal telegraph equation has been relatively slow, with major breakthroughs primarily focusing on the interpretation of its fundamental solutions in terms of stochastic phenomena. The main obstacle is that although the generalized telegraph equation can be transformed into a Volterra type integral equation, the convolution kernel of the latter is not necessarily completely positive. For instance, consider the following time-space nonlocal telegraph equation
	\begin{align}\label{1.2}
		\partial_{t}^{2}(a_1 \ast a_1 \ast (w(\cdot,x)-u(x)))(t) +\gamma\partial_{t}\left(a_1 \ast (w(\cdot,x)-u(x)) \right)(t)+(-\Delta)^{\frac{\beta}{2} } w=g,
	\end{align}
	where $a_1 \in (\mathcal{PC})$. As shown in \cite{J. C. Pozo 20}, even after transforming \eqref{1.2} into the form
	\begin{align}\label{1.3}
		w +b \ast r_\gamma \ast  (- \Delta)^{\frac{\beta}{2}} w= f,
	\end{align}
	numerous counterexamples demonstrate that the kernel $b \ast r_\gamma$ (where $r_\gamma$ is defined in \eqref{2.4}) fails to be completely positive. This critical limitation prevents the direct application of the subordination principle \cite[Chapter 4]{J. Pruss} to derive the relevant solution operators. This issue requires urgent resolution.
	
	Previous studies have primarily focused on normal diffusion (Laplacian operator) and $(\mathcal{PC})$ class kernels. In this work,
we shift our focus to the telegraph equation with a fractional Laplacian and the more general $(\mathcal{PC}^\ast)$ class kernels, which brings together multiple nonlocal phenomena, such as spatial interactions governed by a fractional Laplacian and temporal dynamics with wave-like or memory effects, making it a rich object for theoretical investigation.
In particular, this equation represents an intersection between hyperbolic behavior, anomalous diffusion, and potentially nonlinear effects. The combination of these three elements gives rise to new challenges and behaviors that are not present in classical local models.
Understanding the properties of solutions such as existence, uniqueness, regularity, and long-time behavior remains comparatively underdeveloped, especially when nonlinearities are present. The inclusion of nonlinear terms introduces a delicate competition between three main mechanisms:
\begin{itemize}
  \item the nonlocal diffusion governed by the fractional Laplacian,
  \item the source terms, and
  \item the memory or damping effects, particularly in generalized telegraph-type dynamics.
\end{itemize}
Motivated by the above, we investigate the following fully nonlocal telegraph equation
\begin{align}\label{1.4}
		\begin{cases}
			\partial_t^2 \left(\mathrm{d}a \ast \mathrm{d}a \ast (w-u)\right)(t) +\gamma \partial_t \left(\mathrm{d}a \ast (w-u) \right)(t) +(- \Delta)^{\frac{\beta}{2}}w =- w|w|^{\kappa-1},\\
			w_{t=0}=u(x), \ \partial_t w_{t=0}=0.
		\end{cases}
	\end{align}
	Here, $t>0$, $x \in \mathbb R^3$, $\gamma >0$, $\kappa>1$ and $(-\Delta)^{\frac{\beta}{2}}$ with $\beta \in (1,2]$ is the fractional Laplacian operator. The kernel $a$ is a creep function (see \cite[Definnition 4.4]{J. Pruss}) such that
	there exists a function $b \in L_{loc}^1(\mathbb R_+)$ satisfying
	$(da \ast b)(t) = 1$ for $t>0$, that is, $a$ satisfies the following condition
	\begin{enumerate}
		\item[$(\mathcal{PC}^{\ast})$] The function $a$ admits the representation
		$a(t) = a_0 + \int_0^t a_1(\zeta) \, \mathrm{d}\zeta$ for $t >0$, where $a_0 \geq 0$, $a_1 \in L^1_{loc}(\mathbb R_+)$ is nonnegative and nonincreasing, and there exists $b \in L_{\text{loc}}^1(\mathbb{R}_+)$ such that
		\begin{align}\label{1.5}
			a_0 b(t) + (a_1 \ast b)(t) = 1, \quad t>0.
		\end{align}
	\end{enumerate}
	In this case, we denote by $a \in (\mathcal {PC}^{\ast})$ or $(a,b) \in (\mathcal{PC}^{\ast})$. Obviously, the condition $a_1\in(\mathcal {PC})$ denotes the special case where $a\in (\mathcal{PC}^{\ast })$ with $a_0=0$.

	Similar to the case for \eqref{1.2}, the convolution kernel in the Volterra integral equation obtained from \eqref{1.4} may not be completely positive. To overcome this obstacle, we develop a novel construction of completely positive kernels (see Appendix \ref{app b}), which provides an essential tool for handling nonlocal evolution equations like problem \eqref{1.4} that may fail to be transformed into Volterra integral equations with completely positive kernels. Our analysis reveals that the solution to problem \eqref{1.4} admits a novel representation as
	\begin{align*}
		w(t) = \mathcal C(t)u -\int_0^t \mathcal S(t-\zeta)w(\zeta)|w(\zeta)|^{\kappa-1} \, \mathrm{d}\tau, \quad t \geq 0.
	\end{align*}
	where $\mathcal C(t)$ and $\mathcal S(t)$ are given respectively by
	\begin{align} \label{1.6}
		\mathcal C(t) = - \int_0^\infty  \cos(\theta(-\Delta)^{\frac{\beta}{4} }) \varpi(t,\mathrm{d}\theta), \ \mathcal S(t) = \int_0^\infty  \frac{\sin(\theta(-\Delta)^{\frac{\beta}{4} })}{(-\Delta)^{\frac{\beta}{4} } } \vartheta(t,\mathrm{d}\theta),
	\end{align}
	where $-\varpi(t,\mathrm{d}\theta)$ and $\vartheta(t,\mathrm{d}\theta)$ are positive finite measures induced by the functions $-\varpi(t,\theta)$ and $\vartheta(t,\theta)$, respectively, and in particular, the measure $-w(t, \mathrm{d}\theta)$ satisfies $\int_0^\infty \varpi(t,\mathrm{d}\theta) =1$. Such expressions prove extremely useful for deriving relevant framework space estimates of subordinate wave operators $\mathcal{C}(t)$ and $\mathcal{S}(t)$, as it allows the problem to be reduced to the analysis of the two functions
	\begin{align}\label{1.7}
		c_\delta(t) =\int_0^\infty \theta^{-\delta} \varpi(t,\mathrm{d}\theta) \text{ and } d_\delta(t) = \int_0^\infty \theta^{-\delta} \vartheta(t,\mathrm{d}\theta) \text{ with }  \delta \in (0,1).
	\end{align}
	Under the $(\mathcal{PC})$ condition, it was established in \cite{J. C. Pozo 24} that $c_\delta(t)$ and $d_\delta(t)$ can be expressed as certain integrals involving the relaxation function, and their precise asymptotic behavior was derived. A key observation in our work is that \eqref{1.7} persist when the $(\mathcal{PC})$ condition is relaxed to $(\mathcal{PC}^\ast)$ condition. This insight enables the direct application of existing asymptotic results. Another significant contribution of our work is the first systematic investigation of the almost sure local existence of solutions to problem \eqref{1.4} with regard to random initial conditions. A critical phenomenon is observed in the temporal regularity of the solution with respect to the regularity index $s$ of random initial data $u^\omega \in L^2(\Omega, H^{s,p})$ with $p \in (1,2)$. Specifically, for the subcritical case where $s \in [\frac{2(\kappa-1)}{\kappa+1}, \frac{3(\kappa-1)}{(\kappa+1)})$, the problem \eqref{6.1} admits a unique solution that exists in the weaker space $L^q_tL_x^{1+\kappa}$ with some $q \in (1,\infty)$ almost surely for randomized initial data. However, for the critical and supercritical cases where $s \geq  \frac{3(\kappa-1)}{\kappa+1}$, a unique solution almost surely exists in the stronger space $L^\infty_tL_x^{1+\kappa}$.
	
	Our paper is organized as follows. In the next section, we introduce some notations and function spaces, some deterministic and probabilistic results, and the relaxation function theory. In Section 3, we present an explicit representation of the mild solution to the problem \eqref{1.4}. In section 4, we further establish $L^q_tL^{p'}_x$ estimates for the solution operators. In section 5, using probabilistic methods, we derive a priori estimates for the randomized initial conditions and investigate the resulting averaging phenomena in the context of free evolution. In section 6, we establish probabilistic estimates for the local existence and uniqueness to the considered problem, along with several illustrative examples that support our main results. Appendix \ref{app a} presents the fundamental lemma for deriving the novel solution operators, along with two crucial integrals essential for establishing the operators estimates. Appendix \ref{app b} develops a novel construction method for completely positive kernels. Appendix \ref{app c} establishes fundamental estimates for fractional wave operators.
	
	\section{Preliminaries}
	\subsection{Notations and function spaces}
	Let $f_1 \in \mathcal{S}'(\mathbb R^3)$. We denote it's Fourier transform by either $\tilde {f_1}$ or $\mathcal F(f_1)$, and it's inverse Fourier transform by $\mathcal F^{-1}(f_1)$. Let $\lambda \in \mathbb C$ with $Re(\lambda) >0$. The Laplace transform of a function $f_2:(0,\infty) \to \mathbb C$ or a measure $\mathrm{d}f_3$ is defined as
	$$
	\widehat{f_2}(\lambda)= \mathcal{L}(f_2)(\lambda) = \int_0^\infty f_2(t)e^{-\lambda t}\, \mathrm{d}t, \text{ ~and~ } \mathcal{L}(\mathrm{d}f_3)(\lambda) = \int_0^\infty e^{-\lambda t}\, \mathrm{d}f_3(t),
	$$
	respectively. We denote the inverse Laplace transform by $\mathcal{L}^{-1}$.
	
	Take $s \in \mathbb R$ and $1\leq p,q \leq \infty$. Denote $H^{s,p}(\mathbb R^3)$ as the general fractional Sobolev spaces. Now we introduce the Littlewood–Paley dyadic decomposition. Let $\varrho$ be a nonnegative $C^\infty$ function on $\mathbb{R}^3$ with
	$\text{supp}(\varrho) \subset \{\xi \in \mathbb{R}^3 : |\xi| \le 2\}$ and $\varrho(\xi) = 1$ if $|\xi| \leq 1$. For $j \in \mathbb{N}_+$, let
	$\varrho_j(\xi) = \varrho(2^{-j}\xi) - \varrho(2^{-j + 1}\xi)$ for $\xi \in \mathbb{R}^3$. Then we have $\text{supp}(\varrho_j) \subset \{\xi \in \mathbb{R}^3 : 2^{j-1}\leq |\xi| \leq 2^{j+1}\}$ for $j \in \mathbb N_+$.
	Let $\varrho_0 = \varrho$. The Besov space $B^s_{p,q}(\mathbb R^3)$ is defined as the closure of $C_c^\infty$ in the norm
	\begin{align*}
		\|u\|_{B^s_{p,q}}=\Big(\sum_{j=0}^\infty 2^{qsj}\|\mathcal F^{-1} (\varrho_j \tilde{u})\|^q_{L^p(\mathbb R^3)} \Big)^{\frac{1}{q}}.
	\end{align*}
	As shown in \cite[Theorem 6.4.3]{J. Bergh}, for each $p \in (1,\infty)$ and each $u \in L^p(\mathbb R^3)$, it holds that
	$\|u\|_{L^{p}} \sim \Big\|\Big(\sum\limits_{j=0}^\infty |\mathcal F^{-1} (\varrho_j \tilde{u})|^2 \Big)^{\frac{1}{2}}\Big\|_{L^p(\mathbb R^3)}$. We immediately obtain that
	\begin{align*}
		\|u\|_{H^{s,p}} \sim \Big\|\Big(\sum\limits_{j=0}^\infty |\mathcal F^{-1} (\varrho_j(\xi) (1+|\xi|^2)^{\frac{s}{2} }\tilde{u}(\xi))|^2 \Big)^{\frac{1}{2}}\Big\|_{L^p(\mathbb R^3)} \text{ ~for~ } u \in H^{s,p}(\mathbb R^3).
	\end{align*}
	
	Next, we shall randomize functions by the Littlewood–Paley dyadic decomposition. Let $\{X_j\}_{j \in \mathbb N}$ denote a sequence of mutually independent, zero-mean, real-valued random variables defined on the probability measure space
	$(\Omega, \mathcal A, \mathbb P)$ with their distribution functions denoted by
	$\mu_j$. Assume that there exists $\iota>0$ such that
	\begin{align}\label{2.1}
		\left|\int_{-\infty}^{\infty} e^{\epsilon x}\, \mathrm{d}\mu_j(x)\right| \leq e^{\iota \epsilon^2} ,\quad  \epsilon \in \mathbb R, \ j \in \mathbb N.
	\end{align}
	For a given $u \in H^{s,p}(\mathbb R^3)$ with $p \in (1,2)$, we define its randomization by
	\begin{align}\label{2.2}
		u^\omega = \sum_{j=0}^{\infty} X_j(\omega)\mathcal F^{-1} (\varrho_j \tilde{u}),
	\end{align}
	where the expression above represents the limit of the sequence $\{\sum\limits_{j=0}^N X_j(\omega)\mathcal F^{-1} (\varrho_j \tilde{u})\}$ in $L^2(\Omega, H^{s,p}(\mathbb R^3))$.
	
	Lastly, we outline some crucial functions.
	A function $f \in C^{\infty}\left((0, \infty); \mathbb R \right)$ is termed a completely monotonic function, if $(-1)^jf^{(j)}(\lambda) \geq 0$
	for all $\lambda >0$ and $j \in \mathbb N$. We denote by $\mathcal{CM}$ the class of all such functions. A function
	$f \in C^{\infty}\left((0, \infty); \mathbb R \right)$ is called a Bernstein function if $f(\lambda) \geq 0$ for $\lambda >0$ and $f'(\lambda) \in \mathcal{CM}$. We denote such functions by $\mathcal {BF}$. The profound theory of Bernstein functions can be found in \cite{R. L. Schilling}.
	
	\subsection{Deterministic and probabilistic Preliminaries}
	The following estimate can be found in \cite[Lemma 4]{P. Brenner}. It's proof was originally established by Littman in \cite{W. Littman}.
	\begin{lemma}\label{lemma 2.1}
		Let $Q$ be real, $C^\infty$ in a neighborhood of the support of $f \in C_0^\infty$. Assume that there exist a integer $\varphi \geq 1$ such that the Hessian matrix $H_Q$ satisfies
		\begin{align*}
			\text{rank}(H_Q) \geq \varphi \text{ for all } \xi \in \text{supp}(f).
		\end{align*}
		Then for some integer \( M \), it holds that
		\begin{align*}
			\|\mathcal{F}^{-1}(e^{itQ} f)\|_\infty \leqq C(1 + |t|)^{-\frac{1}{2}\varphi} \sum_{|\alpha| \leqq M} \|D^\alpha f\|_1,
		\end{align*}
		where the constant $C$ depends on the derivatives bounds of $Q$ over $\text{supp}(v)$ and a positive lower bound for the maximal absolute value of $\varphi-th$ order minors of the Hessian matrix $H_Q$ throughout $\text{supp}(v)$.
	\end{lemma}	
	
	We present a criterion for identifying strong
	$(p,q)$ type operators from \cite[Theorem 1.11]{Hormander}.
	\begin{lemma}\label{lemma 2.2}
		 Let $f$ be a measurable function on $\mathbb R^3$ satisfying
		 \begin{align*}
		 	|\{\xi ; \vert f(\xi) \vert \geq \alpha \}| \leq C \alpha^{-\psi} \text{ ~for some constant~ } C \text{ with } 1<\psi<\infty.
		 \end{align*}
		 Then $u \mapsto \mathcal F^{-1} (f(\xi)\tilde u(\xi) )$ is a strong $(p_1,q_1)$ type operator provided
		 \begin{align*}
		 	1<p_1\leq 2 \leq q_1 < \infty,\ \frac{1}{p_1} - \frac{1}{q_1}= \frac{1}{\psi}.
		 \end{align*}
	\end{lemma}
		
	Below, we recall a large deviation estimate from \cite[Lemma 3.1]{N. Burq}.
	\begin{lemma}\label{lemma 2.3}
		Let $\{X_j\}_{j \in \mathbb N}$ denote a sequence of mutually independent, zero-mean, real-valued random variables defined on the probability measure space
		$(\Omega, \mathcal A, \mathbb P)$ with their distribution functions denoted by
		$\mu_j$. Assume that \eqref{2.1} holds.
		Then for every $p_2\geq 2$ and every $\{e_j\} \in l^2(\mathbb N;\mathbb C)$, there exists $L>0$
		such that
		\begin{align*}
			\Big\|\sum_{j=0}^{\infty} X_j(\omega)e_j \Big\|_{L_\omega^{p_2}(\Omega)} \leq L\sqrt{p_2}\Big(\sum_{j=0}^\infty|e_j|^2 \Big)^{\frac{1}{2}}.
		\end{align*}
	\end{lemma}
	
	\subsection{Volterra integral equations}
	Let us recall the Volterra integral equations theory. Let $\gamma \in \mathbb{C}$. For a function $l \in L_{loc}^1(\mathbb R^+)$, we define the relaxation equations as
	\begin{align}
		&s_\gamma(t) + \gamma(l\ast s_\gamma)(t)= 1, \quad t \geq 0. \label{2.3} \\
		&r_\gamma(t) + \gamma(l\ast r_\gamma)(t)= l(t), \quad t>0, \label{2.4}
	\end{align}
	where $s_\gamma, r_\gamma \in L_{loc}^1(\mathbb R^+)$ are the unique solution to the above, respectively, see \cite[Chapter 2, Theorem 3.1]{G. Gripenberg}.
	Sometimes, $s_\gamma$ is also referred to as the scalar resolvent function, while $r_\gamma$ is called the integrated scalar resolvent function. Let $f \in L_{loc}^p(\mathbb R_+)$ for some $p\geq 1$, then the unique solution to the following Volterra integral equation
	\begin{align}\label{2.5}
		v_\gamma(t) + \gamma(l\ast v_\gamma)(t) = f(t), \quad t>0.
	\end{align}
	is given by
	\begin{align}\label{2.6}
		v_\gamma(t) = f(t) - \gamma (r_\gamma \ast f)(t), \quad t>0,
	\end{align}
	and $v_\gamma \in L_{loc}^p(\mathbb R_+)$, see \cite[Chapter 2, Theorem 3.5]{G. Gripenberg}. Furthermore, if $f \in C(\mathbb R_+)$, then $v_\gamma \in C(\mathbb R_+)$ and the previous form of $v_\gamma(t)$ also holds for $t=0$, see \cite[Chapter 2, Theorem 3.5]{G. Gripenberg}.
	
	Generally, if $s_\gamma$ and $r_\gamma$ is nonnegative for all $\gamma \geq 0$, then
	we say that $l$ is a completely positive function, see \cite[Definition 1.1]{P. Clement 81}. As shown in \cite[Theorem 2.2]{P. Clement 81}, complete positivity admits the following equivalent condition: there exist a nonnegative nonincreasing function $k_1 \in L_{loc}^1(\mathbb R_+)$ and $k_0\geq 0$ such that
	\begin{align}\label{2.7}
		k_0l(t) + (k_1\ast l)(t) = 1, \quad t>0.
	\end{align}
	Thus, $(a,b) \in (\mathcal{PC}^{\ast})$ implies that $b$ is completely positive.
	It should be mentioned that a completely positive function must be nonnegative, see  \cite[Propposition 2.1(i)]{P. Clement 81}.
	Further equivalent conditions regarding the complete positivity of $b$ is presented in \cite[Proposition 4.5]{J. Pruss}.
	
	Presented below are several key properties of $s_\gamma$ and $r_\gamma$.
	\begin{proposition}\label{prop 2.1}
		Assume that $l$ is completely positive and $\gamma \geq 0$. Then the following hold:
		\begin{itemize}
			\item[(i)] The function $s_\gamma(t)$ has the following representation
			\begin{align*}
				s_\gamma(t) = 1 - \gamma \int_0^t r_\gamma(\zeta)\, \mathrm{d}\zeta = k_0 r_\gamma(t)+ (k_1 \ast r_\gamma)(t)=(\mathrm{d}k \ast r_\gamma)(t), \quad t \geq 0,
			\end{align*}
			which implies that $s_\gamma$ is continuous on $\mathbb R_+$ and differentiable on $(0,\infty)$.
			\item[(ii)] For $t>0$, it holds
			\begin{align*}
				\frac{1+k_0\gamma [k_1(t)]^{-1}r_\gamma(t)}{1+\gamma [k_1(t)]^{-1}} \leq s_\gamma(t) \leq \frac{1}{1+\gamma (1\ast l)(t)}.
			\end{align*}
			\item[(iii)] For $t>0$, it holds
			\begin{align*}
				0 \leq r_\gamma(t) \leq l(t), \ (1\ast r_\gamma)(t) \leq \frac{(1\ast l)(t)}{1+\gamma (1\ast l)(t)}.
			\end{align*}
			 Moreover if $l$ is nonincreasing, then
			\begin{align*}
				r_\gamma(t) \leq \frac{ l(t)}{1+\gamma(1\ast l)(t)}.
			\end{align*}
		\end{itemize}
		\begin{proof}
			The justification for (i) and (ii) can be founded in \cite[Proposition 2.1]{X. Huang}. For (iii), the proof is completely analogous to that in \cite[Lemma 5.4]{J.C. Pozo}, so we omit it.
		\end{proof}
	\end{proposition}
	
	\section{Novel solution representation}
	Assume that $(a,b)\in (\mathcal {PC}^{\ast})$. Without loss of rigor, we can specifically choose $l=b$ in both \eqref{2.4} and \eqref{2.5}. Then, applying Lemma \ref{lemma B.1} and \eqref{B.3} yields that $h_{\frac{1}{2},\frac{1}{2},\gamma} \in (\mathcal{PC}^{\ast})$ whose Laplace transform is given by
	\begin{align*}
		\widehat{h}_{\frac{1}{2},\frac{1}{2},\gamma}(\lambda):=\pi \sqrt{\widehat{b}(\lambda)\widehat{r}_\gamma(\lambda)},
	\end{align*}
	and there exists a nonnegative nonincreasing function $m_1$ and $m_0\geq 0$ such that
	\begin{align}\label{3.1}
		m_0h_{\frac{1}{2}, \frac{1}{2},\gamma}+(m_1 \ast h_{\frac{1}{2},\frac{1}{2},\gamma})(t) =1, \quad t>0,
	\end{align}
	For simplicity, let $h_\gamma= {\pi}^{-\frac{1}{2} } h_{\frac{1}{2},\frac{1}{2},\gamma}$. Then $h_\gamma$ is completely positive and
	\begin{align}\label{3.2}
		(h_\gamma \ast h_\gamma)= (b \ast r_\gamma)(t), \quad t>0.
	\end{align}
	For each $\gamma \geq 0$, consider the following relaxation equations:
	\begin{align}
		&z(t,\nu) + \nu(h_\gamma \ast z(\cdot,\nu))(t)= 1, \quad t \geq 0, \label{3.3}\\
		&n_\nu(t) + \nu(h_\gamma \ast n(\cdot,\nu))(t)= h_\gamma(t), \quad t>0, \label{3.4},
	\end{align}
	where $\nu \in \mathbb C$. By Lemma \ref{lemma A.1}, there exist unique positive finite measures $-\varpi(t,\mathrm{d}\theta)$ and $\vartheta(t,\mathrm{d}\theta)$ such that for all $\nu \in \{z \in \mathbb {C}: Re(z)\geq -a_\infty \} $, the functions $z_\nu(t)$ and $ n_\nu(t)$ admit the representations:
	\begin{align}\label{3.5}
		z(t,\nu) = -\int_0^\infty e^{-\nu \theta} \varpi(t,\mathrm{d}\theta),\ n(t,\nu) = \int_0^\infty e^{-\nu \theta} \vartheta(t,\mathrm{d}\theta),  \quad t>0.
	\end{align}
	
	Consider the initial value problem for the generalized telegraph equation:
	\begin{align*}
		\begin{cases}
			\partial_t^2 \left(\mathrm{d}a \ast \mathrm{d}a \ast (w-u)\right)(t) +\gamma \partial_t (\mathrm{d}a \ast (w-u))(t) +(- \Delta)^{\frac{\beta}{2}}w = g(t,x), \quad t>0,\\
			w_{t=0}=u(x), \ \partial_t w_{t=0}=0, \quad x \in \mathbb R^3.
		\end{cases}
	\end{align*}
	
	Convolving both sides of the equation with
	$ b \ast r_\gamma$ yields
	\begin{align*}
		&\quad r_\gamma \ast b\ast\partial_t^2 \left(\mathrm{d}a \ast \mathrm{d}a \ast (w-u)\right) +\gamma r_\gamma \ast b\ast \partial_t (\mathrm{d}a \ast (w-u))\\
		&= r_\gamma \ast \partial_t (b\ast \partial_t \mathrm{d}a \ast \mathrm{d}a\ast (w-u)) + \gamma r_\gamma \ast \partial_t (b \ast \mathrm{d}a \ast (w-u))\\
		&= r_\gamma \ast  \partial_t (\mathrm{d}a\ast (w-u)) + \gamma r_\gamma \ast \partial_t (1  \ast (w-u))\\
		&=r_\gamma \ast \partial_t((\mathrm{d}a+\gamma 1)\ast (w-u))\\
		&=\partial_t (r_\gamma \ast (\mathrm{d}a+\gamma 1)\ast (w-u))\\
		&=\partial_t (1\ast(w-u))\\
		&=w-u.
	\end{align*}
	Here, we employ the fact that
	$((\mathrm{d}a +\gamma 1)\ast r_\gamma)(t) = 1$ for $t\geq 0$ and $\gamma \in \mathbb R_+$ (from Proposition \ref{prop 2.1}(i)). It directly derives that
	\begin{align*}
		w(t) +(b \ast r_\gamma \ast (-\Delta)^{\frac{\beta}{2}}w)(t) = u + (b \ast r_\gamma \ast g)(t), \ w_{t=0} =u(x), \quad  t>0.
	\end{align*}
	Taking the Fourier transform with respect to the spatial variable $x$ on both sides, we obtain
	\begin{align}\label{3.6}
		\begin{cases}
			\tilde{w}(t,\xi) + |\xi|^\beta(b \ast r_\gamma \ast \tilde{w}(\cdot,\xi)) =\tilde{w}_0(\xi) + (b \ast r_\gamma \ast \tilde{g}(\cdot,\xi))(t),\quad \ t>0,\\
			\tilde {w}_{t=0} =\tilde{u}(\xi),
		\end{cases}
	\end{align}
	
	We now define two functions $C$ and $S$ as follows:
	\begin{align}
		&C(t,\xi):=
			Re(z(t,i|\xi|^{\frac{\beta}{2}}))=- \int_0^\infty \cos(\theta|\xi|^{\frac{\beta}{2}}) \varpi(t,\mathrm{d}\theta), \quad t\geq 0, \ \xi \in \mathbb R^3,
		 \label{3.7}\\
		&S(t,\xi):=\begin{cases}
			\frac{Im(n(t,-i|\xi|^{\frac{\beta}{2}}))}{|\xi|^{\frac{\beta}{2}}} =\int_0^\infty |\xi|^{-\frac{\beta}{2}} \sin(\theta|\xi|^{\frac{\beta}{2}}) \vartheta(t,\mathrm{d}\theta), \quad t>0, \ \xi \in \mathbb R^3\setminus\{0\},\\
			(b \ast r_\gamma)(t), \quad t>0,\ \xi =0.
		\end{cases} \label{3.8}
	\end{align}
	Setting $\nu=i|\xi|^{\frac{\beta}{2} }$ into \eqref{3.3}, we obtain
	\begin{align*}
		z(t,i|\xi|^{\frac{\beta}{2}}) + i|\xi|^{\frac{\beta}{2}} (h_\gamma \ast z(\cdot,i|\xi|^{\frac{\beta}{2}})) = 1, \quad t\geq 0,
	\end{align*}
	This equation equivalently writes as
	\begin{empheq}[left=\empheqlbrace]{align}
		&Re(z(t,i|\xi|^{\frac{\beta}{2}})) - |\xi|^{\frac{\beta}{2}} (h_\gamma \ast Im(z(\cdot,i|\xi|^{\frac{\beta}{2}})))(t) = 1, \quad t\geq 0, \label{3.9}\\
		&Im(z(\cdot,i|\xi|^{\frac{\beta}{2}}))+ |\xi|^{\frac{\beta}{2}} (h_\gamma \ast Re(z(\cdot,i|\xi|^{\frac{\beta}{2}})))(t)=0,  \quad t \geq 0.\label{3.10}
	\end{empheq}
	Substituting \eqref{3.10} into \eqref{3.9}, then by \eqref{3.2} and \eqref{3.7} we have
	\begin{align}\label{3.11}
		C(t,\xi) + |\xi|^\beta(b \ast r_\gamma \ast C(\cdot,\xi))(t) =1, \quad t\geq 0.
	\end{align}
	Similarly, The following equation
	\begin{align*}
		n(t,-i|\xi|^{\frac{\beta}{2}}) - i|\xi|^{\frac{\beta}{2}} (h_\gamma \ast n(\cdot,-i|\xi|^{\frac{\beta}{2}})) = h_\gamma(t), \quad t> 0.
	\end{align*}
	is equivalent to
	\begin{empheq}[left=\empheqlbrace]{align*}
		&Re(n(t,-i|\xi|^{\frac{\beta}{2}})) + |\xi|^{\frac{\beta}{2}} (h_\gamma \ast Im(n(\cdot,-i|\xi|^{\frac{\beta}{2}})))(t) = h_\gamma(t) , \\
		&Im(n(t,-i|\xi|^{\frac{\beta}{2}}))- |\xi|^{\frac{\beta}{2}} (h_\gamma \ast Re(n(\cdot,-i|\xi|^{\frac{\beta}{2}})))(t)=0, \quad t > 0.
	\end{empheq}
	By analogy with the derivation of \eqref{3.10}, one can derive that
	\begin{align}\label{3.12}
		S(t,\xi) + |\xi|^\beta (b \ast r_\gamma \ast S(\cdot,\xi))(t) = (b \ast r_\gamma)(t), \quad t> 0.
	\end{align}
	Finally, from \eqref{2.6}, \eqref{3.6}, \eqref{3.11} and \eqref{3.12} we obtain
	\begin{align*}
		\tilde{w}(t,\xi) &= \tilde{u}(\xi)\left(1  - |\xi|^\beta (1\ast S(\cdot,\xi))(t) \right)\\
		&\quad +(b \ast r_\gamma - |\xi|^\beta b \ast r_\gamma \ast S(\cdot,\xi)) \ast \tilde{g}(\cdot,\xi)(t)\\
		&=C(t,\xi)\tilde{u}(\xi) + \int_0^t S(t-\zeta,\xi) \tilde{g}(\zeta,\xi) \, \mathrm{d}\zeta.
	\end{align*}
	then the uniqueness of Fourier transform leads to
	\begin{align*}
		w(t) = {\mathcal F}^{-1}(C(t,\cdot) \tilde{u} ) + \int_0^t {\mathcal F}^{-1}\left(S(t-\zeta, \cdot)\tilde{g}(\zeta,\cdot)\right)\, \mathrm{d}\zeta, \quad t>0.
	\end{align*}
	Let us define the following operators:
	\begin{align*}
		\begin{aligned}
			&\mathcal C(t)u := \mathcal F^{-1} (C(t,\xi)u(\xi))(x)= - \int_0^\infty  \cos(\theta(-\Delta)^{\frac{\beta}{4} }) \varpi(t,\mathrm{d}\theta), \quad t\geq 0,\\
			&\mathcal S(t)u :=\mathcal F^{-1} (S(t,\xi)u(\xi))(x) =  \int_0^\infty  \frac{\sin(\theta(-\Delta)^{\frac{\beta}{4} })}{(-\Delta)^{\frac{\beta}{4} } } \vartheta(t,\mathrm{d}\theta),\quad t>0.
		\end{aligned}
	\end{align*}
	Consequently, the solution to the aforementioned problem can be expressed as
	\begin{align*}
		w(t) = \mathcal C(t) u + \int_0^t \mathcal S(t-\zeta, \cdot)g(\zeta)\, \mathrm{d}\zeta, \quad t\geq 0.
	\end{align*}
	
	\section{Space-time estimates for the fully nonlocal telegraph equations}
	We begin by establishing the space-time estimates for the operator $\mathcal C(t)$. For this purpose, we impose the following assumptions on $\beta, s, p$.
	\begin{enumerate}
		\item[$(\mathcal{H}_1)$] $s\in (0,3)$, $p \in (1, \frac{6}{3+s}] \cap [\frac{4}{2+s}, \frac{6}{3+s}]$ if $\beta =2$, and $s \in (0, 2-\frac{\beta}{2})\cup[3-\frac{\beta}{2} ,3)$, $p \in (1, \frac{6}{3+s}] \cap [ 2\frac{12-3\beta}{12-3\beta +4s}, \frac{6}{3+s}]$ if $\beta \in (1,2)$.
	\end{enumerate}
	
	For $p\in (1,2)$, let $p'$ denote the conjugate exponent of $p$. For convenience, denote $\tau_1(\beta,p,s)= \frac{6}{\beta}(\frac{1}{p} - \frac{1}{p'})-\frac{2s}{\beta}$. Direct computation confirms that under the hypothesis $({\mathcal H}_1)$ $\tau_1(\beta,p,s)$ takes values in $[0,1)$ . We put
	\begin{align*}
		\rho_{1,\delta}(t)= -\int_0^\infty \theta^{-\delta} \varpi(t,\mathrm{d}\theta), \quad t>0 , \ \delta \in [0,1).
	\end{align*}
	
	\begin{lemma}\label{lemma 4.1}
		Assume that $(a,b)\in (\mathcal{PC}^{\ast})$ and $T>0$. If $(\mathcal{H}_1)$ holds and $0\leq q \tau_1(\beta,p,s) < 1$ for some $q \in [1, \infty]$. Then there exists a continuous, nondecreasing function $M_1:[0,\infty) \to [0,\infty)$ such that for all $u \in H^{s,p}(\mathbb R^3)$, the following estimate holds:
			\begin{align}\label{4.1}
				\|\mathcal C(t)u\|_{L^q([0,T]; L^{p'}(\mathbb R^3))} \leq M_1(T) \|u\|_{H^{s,p}}.
			\end{align}
		\begin{proof}
			Note that $\beta,s,p$ now satisfy the condition of Proposition \ref{prop C.1}. An application of \eqref{3.7} with Proposition \ref{prop C.1} and Lemma \ref{lemma A.2} yields
			\begin{align*}
				\|\mathcal C(t)u \|_{L^{p'}(\mathbb R^3)}
				&=\left\|\mathcal F^{-1}\left(C(t,\xi)\tilde u(\xi) \right) \right\|_{L^{p'}} \\
				&=\left\|-\int_0^\infty \mathcal F^{-1} \left(\cos( \theta|\xi|^{\frac{\beta}{2}}) \tilde{u}(\xi) \right) \varpi(t,\mathrm{d}\theta)  \right\|_{L^{p'}}  \\
				&\leq-\int_0^\infty\left\|\mathcal F^{-1} \left(\cos( \theta|\xi|^{\frac{\beta}{2}}) \tilde{u}(\xi) \right) \right\|_{L^{p'}} \varpi(t,\mathrm{d}\theta)\\
				&\lesssim -\int_0^\infty \theta^{\frac{2s}{\beta} - \frac{6}{\beta}(\frac{1}{p} - \frac{1}{p'}) } \varpi(t,\mathrm{d}\theta) \|u\|_{H^{s,p}} \\
				&\lesssim \rho_{1,\tau_1(\beta,p,s)}(t)
				\|u\|_{H^{s,p}} \text{ ~for~ } t > 0.
			\end{align*}
			If $\tau_1(\beta,p,s) =0$, which occurs when  $s=3(\frac{1}{p}-\frac{1}{p'})$, then Lemma \ref{A.1}(i) implies that  $\rho_{\tau_1(\beta,p,s)}(t)|_{\tau_1(\beta,p,s)=0} = 1$ for $t \geq 0$. For the case $\tau_1(\beta,p,s)\in (0,1)$, Lemma \ref{A.2} directly yields that $\rho_{\tau_1(\beta,p,s)}(t)\lesssim [(1 \ast h_\gamma)(t)]^{-\tau_1(\beta,p,s) }$ for $t>0$. Therefore, we conclude that
			\begin{align}\label{4.2}
				\|\mathcal C(t)u \|_{L^{p'}(\mathbb R^3)}
				\lesssim [1\ast h_\gamma(t)]^{-\tau_1(\beta,p,s)} \|u\|_{H^{s,p}} \text{ ~for~ } t > 0,
			\end{align}
			provided $(\mathcal{H}_1)$ holds. It follows from \eqref{3.1} that
			$[(1\ast h_\gamma)(t)]^{-1}\lesssim t^{-1}[m_0+ (1\ast m_1)(t)]$ for $t>0$,
			which shows that
			\begin{align}\label{4.3}
				\|\mathcal C(t)u \|_{L^{p'}}\lesssim [m_0+(1\ast m_1)(t)]^{\tau_1(\beta,p,s) } t^{-\tau_1(\beta,p,s)} \|u\|_{H^{s,p}}, \quad t>0.
			\end{align}
			We define $M_1:\mathbb R_+ \to \mathbb R_+$ as follow:
			$M_1(T)= K_1[m_0+(1\ast m_1)(T)]^{\tau_1(\beta,p,s) }T^{\frac{1}{q}-\tau_1(\beta,p,s)}$,
			where $K_1$ is sufficiently large. This completes the proof of \eqref{4.1}.
		\end{proof}
	\end{lemma}
	
	\begin{remark}\label{rem 4.1}
		\eqref{4.2} reveals a critical index phenomenon, where the critical index is given by $s_{1,crit}= 3(\frac{2}{p}-1)$.
		Under hypothesis $(\mathcal{H}_1)$, $\|\mathcal C(t)u\|_{L^{P'}}$ exhibits a time decay estimate when $s<s_{1,crit}$, while it satisfies a uniformly time bounded estimate when $s=s_{1,crit}$.
	\end{remark}
	
	We now proceed to derive the space-time estimates associated with the operator $\mathcal S(t)$. Our analysis requires the following parameter constraints for $\beta$ and $p$.
	\begin{enumerate}
		\item[$(\mathcal{H}_2)$] $p \in [\frac{4}{3}, \frac{3}{2})$ if $\beta =2$, and $p \in [2\frac{12-3\beta}{12-\beta},\frac{12}{6+\beta})$ if $\beta \in (1,2)$.
	\end{enumerate}
	
	\begin{enumerate}
		\item[$(\mathcal{H}_3)$] $p \in [\frac{3}{2}, 2]$ if $\beta =2$, and $ p \in [\frac{12}{6+\beta},2]$ if $\beta \in (1,2)$.
	\end{enumerate}
	
	For convenience, denote $\tau(\beta,p)= \frac{6}{\beta}(\frac{1}{p} - \frac{1}{p'})-1 $. Direct computation confirms that under the $(\mathcal{H}_2)$ or $(\mathcal{H}_3)$, $\tau(\beta,p)$ takes values in $(0,1)$ and $[-1,0]$, respectively. We put
	\begin{align*}
		&\sigma(t) =  \int_0^\infty \theta^{-\tau(\beta,p)} \vartheta(t,\mathrm{d}\theta), \quad t>0 , \ \tau(\beta,p) \text{ ~satisfies~} (\mathcal{H}_2) \text{ or } (\mathcal{H}_3).\\
		&\rho_\delta(t)= \int_0^\infty \theta^{-\delta} \vartheta(t,\mathrm{d}\theta), \quad t>0 , \ \delta \in [0,1).
	\end{align*}
	
	\begin{lemma}\label{lemma 4.2}
		Let $T>0$, $(a,b)\in (\mathcal{PC}^{\ast})$ and $b$ is nonincreasing. Assume that one of the following holds:
		\begin{itemize}
			\item[(i)] $a_0 = 0$, and one of the following conditions holds:
			\begin{itemize}
				\item[(1)] $\beta, p$ satisfy $(\mathcal{H}_2)$,
				and $b \in L_{loc}^{p_0}(\mathbb R_+)$ for some $p_0>\frac{q_0(1-\tau(\beta,p))}{1-q_0\tau(\beta,p)}$ with some $q_0\geq 1$ satisfying $0<q_0\tau(\beta,p)<1$.
				\item[(2)] $\beta,p$ satisfy $(\mathcal{H}_3)$, and $b \in L_{loc}^{q_0}(\mathbb R_+)$ for some $q_0 \geq 1$.
			\end{itemize}
			\item[(ii)] $a_0>0$, and one of the following conditions holds:
			\begin{itemize}
				\item[(3)]$\beta, p$ satisfy $(\mathcal{H}_2)$ and
				$0< \tau(\beta,p)q_0<1$ for some $q_0\geq 1$.
				\item[(4)]$\beta,p$ satisfy $(\mathcal{H}_3)$ and $q_0 \in[1,\infty]$.
			\end{itemize}
		\end{itemize}
		Then for any exponents $1\leq r\leq q \leq \infty$ satisfying $1+\frac{1}{q}=\frac{1}{q_0}+\frac{1}{r}$, there exists a continuous, nondecreasing function $M:[0,\infty) \to [0,\infty)$ with $M(0)=0$, such that
		\begin{align}\label{4.4}
			\left\|\int_0^t \mathcal S(t-\zeta) v(\zeta)\, \mathrm{d}\zeta \right\|_{L^q([0,T]; L^{p'}(\mathbb R^3))} \leq M(T) \|v\|_{L^r([0,T]; L^{p}(\mathbb R^3))}.
		\end{align}
		\begin{proof}
			Prior to proving these conclusions, we first establish an estimate:
			\begin{align}\label{4.5}
				\rho_{\delta} (t) \lesssim b(\frac{t}{2}) [(1\ast b)(\frac{t}{2})]^{ -\delta}, \quad t>0, \ \delta \in [0,1).
			\end{align}
			Since $h_\gamma$ is completely positive, an application of Lemma \ref{lemma A.2} yields
			\begin{align*}
				\widehat{\rho}_{\delta}(\lambda) &= \Gamma(1- \delta) \widehat{h}_\gamma(\lambda)^{1-\delta}\\
				&=\Gamma(1- \delta) [\widehat{b}(\lambda)]^{\frac{1-\delta}{2}} [\widehat{r}_\gamma(\lambda)]^{\frac{1-\delta}{2}}\\
				&=\frac{\Gamma(1-\delta) }{\Gamma\left(1 -\frac{1+\delta}{2}\right)^2}\Gamma\left(1 -\frac{1+\delta}{2}\right)^2  [\widehat{b}(\lambda)]^{1-\frac{1+\delta}{2}} [\widehat{r}_\gamma(\lambda)]^{1-\frac{1+\delta}{2}}\\
				&=\frac{\Gamma(1-\delta) }{\Gamma\left(1 -\frac{1+\delta}{2}\right)^2} \widehat{h}_{\frac{1+\delta}{2},\frac{1+\delta}{2} ,\gamma} (\lambda) \text{ ~for~ } \delta \in [0,1).
			\end{align*}
			Consequently, we infer that \eqref{4.5} holds by Lemma \ref{lemma B.1}. Hence, it is sufficient to prove that
			\begin{align}\label{4.6}
				\sigma(t) \lesssim b(\frac{t}{4}) \left[(1\ast b)(\frac{t}{4})\right]^{-\tau(\beta,p)} \text{ ~for~ } t >0,
			\end{align}
			provided that $(\mathcal{H}_2)$ or $(\mathcal{H}_3)$ holds. Indeed, \eqref{4.6} is directly obtained from \eqref{4.5} for the case that $\beta, p$ satisfying  $(\mathcal{H}_2)$.
			Given the condition $(\mathcal{H}_3)$, if $\delta(\beta,p) \in [-1,0)$, which corresponds to $p \in ( \frac{3}{2},2]$ when $\beta =2 $ and $p \in (\frac{12}{6+\beta}, 2]$ for $\beta \in (1,2)$, we deduce that
			\begin{align*}
				\widehat{\sigma}(\lambda) &= \int_0^\infty \theta^{-\tau(\beta,p)} \, \widehat{\vartheta}(\lambda,\mathrm{d}\theta) = \int_0^\infty \theta^{-\tau(\beta,p)} e^{-{\theta}[\widehat{h_\gamma}(\lambda)]^{-1}}\, \mathrm{d}\theta\\
				&=\Gamma(1-\tau(\beta,p)) [\widehat{h_\gamma}(\lambda)]^{1- \tau(\beta,p)}\\
				&=\frac{\Gamma(1-\tau(\beta,p))}{\Gamma(-\tau(\beta,p ))}\widehat{h}_\gamma(\lambda)\widehat{\rho}_{1+\tau(\beta,p),\gamma}(\lambda),
			\end{align*}
			which implies $\sigma(t)  =\frac{\Gamma(1-\tau(\beta,p))}{\Gamma(-\tau(\beta,p ))} \left(h_\gamma \ast \rho_{1+\tau(\beta,p)} \right)(t)$.
			Therefore, we immediately obtain
			\begin{align*}
				\sigma(t) &\lesssim (h_\gamma \ast \rho_{1+\tau(\beta,p)})(t)\\
				&\lesssim \int_0^{\frac{t}{2}} b(\frac{t-\zeta}{2}) b(\frac{\zeta}{2}) \left[(1\ast b)(\frac{\zeta}{2}) \right]^{-(1+\tau(\beta,p) )}  \, \mathrm{d}\zeta\\
				&\quad + \int_{\frac{t}{2}}^t b(\frac{t-\zeta}{2}) b(\frac{\zeta}{2}) \left[(1\ast b)(\frac{\zeta}{2}) \right]^{-(1+\tau(\beta,p) )}  \, \mathrm{d}\zeta \\
				&\lesssim b(\frac{t}{4}) \left[(1\ast b)(\frac{t}{4})\right]^{-\tau(\beta,p)}, \quad  t >0.
			\end{align*}
			If $\tau(\beta,p) = 0$, which corresponds to $p = \frac{3}{2}$ when $\beta =2 $ and $p =\frac{12}{6+\beta}$ for $\beta \in (1,2)$, Lemma \ref{lemma A.1}(i) directly yields the estimate
			$\sigma(t) \lesssim b(t)$.
		
			Let us first consider the case (i). According to Proposition \ref{prop C.2} and \eqref{4.6}, it follows that
			\begin{align}\label{4.7}
				\begin{aligned}
					\|\mathcal S(t)u\|_{L^{p'}(\mathbb R^3)} &\leq \int_0^\infty \left\|\mathcal F^{-1}\left(|\xi|^{-\frac{\beta}{2} } \sin t|\xi|^{\frac{\beta}{2} } \tilde{u}(\xi) \right)(\cdot) \right\|_{L^{p'}(\mathbb R^3)} \, \vartheta(t, \mathrm{d}\theta)\\
					&\lesssim  \sigma(t)\|u\|_{L^p(\mathbb R^3) }.
				\end{aligned}
			\end{align}
			Due to $b$ is nonnegative nonincreasing, we derive that
			\begin{align}\label{4.8}
				tb(t) \leq (1 \ast b)(t), \quad t>0.
			\end{align}
			Under the condition (1), since $\tau(\beta,p) \in (0,1]$, one can derive from \eqref{4.6} and \eqref{4.8} that
			\begin{align*}
				\sigma(t)
				\lesssim t^{-\tau(\beta,p)} b(\frac{t}{4})^{1-\tau(\beta,p)}, \quad t>0.
			\end{align*}
			Since $b \in L_{loc}^{p_0}(\mathbb R_+)$ with $p_0 >\frac{q_0(1-\tau(\beta,p))}{1-q_0\tau(\beta,p) }$, applying the generalized H\"{o}lder's inequality ($\frac{1}{q_0}=\frac{1}{\frac{p_0}{p_0-(1-\tau(\beta,p))} } + \frac{1}{\frac{p_0}{1-\tau(\beta,p)} } $) to the above, we can obtain that
			\begin{align}\label{4.9}
				\|\sigma(t)\|_{L^{q_0}([0,T])} \lesssim T^{\frac{1}{q_0}-\frac{1-\tau(\beta,p)}{p_0} - \tau } \|b\|_{L^{p_0}([0,\frac{T}{4}])}.
			\end{align}
			Here, $\frac{1}{q_0}-\frac{1-\tau(\beta,p)}{p_0} - \tau>0$ due to $p_0 >\frac{q_0(1-\tau(\beta,p))}{1-q_0\tau(\beta,p) }$. It follows from \eqref{4.7}, \eqref{4.9} and Young's inequality, we get
			\begin{align*}
				&\quad \left\|\int_0^t \mathcal S(t-\zeta) v(\zeta)\, \mathrm{d}\zeta \right\|_{L^q([0,T]; L^{p'}(\mathbb R^3))}\\
				&\leq
				\left \|\int_0^t\left  \|\mathcal{S}(t-\zeta)v(\zeta) \right\|_{L^{p'}(\mathbb R^3)}\, \mathrm{d}\zeta \right\|_{L^q ([0,T])}\\
				&\lesssim\left \|\int_0^t  \sigma(t-\zeta) \left\|v(\zeta) \right\|_{L^{p}(\mathbb R^3)}\, \mathrm{d}\zeta \right\|_{L^q ([0,T])}\\
				&\lesssim \|\sigma\|_{L^{q_0}([0,T])} \|v\|_{L^{r}([0,T]; L^{p}(\mathbb R^3))}\\
				&\lesssim T^{\frac{1}{q_0}-\frac{1-\tau(\beta,p)}{p_0} - \tau } \|b\|_{L^{p_0}([0,\frac{T}{4}])}\|v\|_{L^{r}([0,T]; L^{p}(\mathbb R^3))}.
			\end{align*}
			Given the condition (2), since $\tau(\beta,p) \in [-1,0]$, one can derive from \eqref{4.6} that
			\begin{align*}
				\|\sigma(t)\|_{L^{q_0}([0,T])} \lesssim [(1\ast b)(\frac{T}{4})]^{-\tau(\beta,p)} \|b\|_{L^{q_0}([0,\frac{T}{4}]) }.
			\end{align*}
			This along with \eqref{4.7} and Young's inequality gives
			\begin{align*}
				\quad \left\|\int_0^t \mathcal S(t-\zeta) v(\zeta)\, \mathrm{d}\zeta \right\|_{L^q([0,T]; L^{p'}(\mathbb R^3))}
				\lesssim \|b\|^{-\tau(\beta,p)}_{L^{1}([0,\frac{T}{4}])} \|b\|_{L^{q_0}([0,\frac{T}{4}])}\|v\|_{L^{r}([0,T]; L^{p}(\mathbb R^3))}.
			\end{align*}
			
			Consider the case (ii). It follows from \eqref{1.5} and $a_0>0$ that $b(t) \in [0, \frac{1}{a_0}]$ for $t\geq 0$ and $[(1\ast b)(t)]^{-1} \leq t^{-1}[a_0+(1\ast a_1)(t)]$ for $t>0$. Immediately, under the condition (3), one can derive from \eqref{4.6} that
			\begin{align*}
				\sigma(t) \lesssim [(1\ast b)(\frac{t}{4})]^{-\tau(\beta,p)} \lesssim \frac{[a_0+(1\ast a_1)(\frac{t}{4})]^{\tau(\beta,p)}}{t^{\tau(\beta,p)}}.
			\end{align*}
			Then following an argument analogous to the proof of Lemma \ref{lemma 4.1}, we establish
			\begin{align*}
				\left\|\int_0^t \mathcal S(t-\zeta) v(\zeta)\, \mathrm{d}\zeta \right\|_{L^q([0,T]; L^{p'}(\mathbb R^3))} \lesssim [a_0+(1\ast a_1)(\frac{T}{4})]^{\tau(\beta,p)}T^{\frac{1}{q_0}-\tau(\beta,p)} \|v\|_{L^r([0,T]; L^{p}(\mathbb R^3))}.
			\end{align*}
			Given the condition (4), from $\tau(\beta,p) \in [-1,0]$ and $b(t) \in [0, \frac{1}{a_0}]$ for $t\geq 0$ we have $\sigma(t)\lesssim t^{-\tau(\beta,p)}$. Then Young's inequality shows that
			\begin{align*}
				\left\|\int_0^t \mathcal S(t-\zeta) v(\zeta)\, \mathrm{d}\zeta \right\|_{L^q([0,T]; L^{p'}(\mathbb R^3))} \lesssim T^{\frac{1}{q_0}-\tau(\beta,p)} \|v\|_{L^r([0,T]; L^{p}(\mathbb R^3))}.
			\end{align*}
			
			Finally, we define $M:\mathbb R_+ \to \mathbb R_+$ as follows:
			\begin{align*}
				M(T)= \begin{cases}
					K_2T^{\frac{1}{q_0}-\frac{1-\tau(\beta,p)}{p_0} - \tau(\beta,p) } \|b\|_{L^{p_0}([0,\frac{T}{4}])}, \text{ ~if (1) holds},\\
					K_2\|b\|_{L^{q_0}([0,\frac{T}{4}])}\|b\|^{-\tau(\beta,p)}_{L^{1}([0,\frac{T}{4}])}, \text{ ~if (2) holds},\\
					K_2[a_0+(1\ast a_1)(\frac{T}{4})]^{\tau(\beta,p)}T^{\frac{1}{q_0}-\tau(\beta,p)} \text{ ~if (3) holds},\\
					K_2T^{\frac{1}{q_0}-\tau(\beta,p)}\text{ ~if (4) holds},
				\end{cases}
			\end{align*}
			where $K_2$ is sufficiently large. This completes the proof.
		\end{proof}
	\end{lemma}
	
	\begin{remark}\label{rem 4.2}
		In fact, under the minimal assumptions that
		$(a,b) \in (\mathcal {PC}^{\ast})$ and $b$ is nonincreasing, we can conclude that
		\begin{align*}
			\left\|\int_0^t \mathcal S(t-\zeta) v(\zeta)\, \mathrm{d}\zeta \right\|_{L^q([0,T]; L^{p'}(\mathbb R^3))} \leq K_2 \|b\|^{1-\tau(\beta,p)}_{L^1([0,\frac{T}{4}])}\|v\|_{L^q([0,T]; L^{p}(\mathbb R^3))} \text{ ~for~ } q \in [1,\infty],
		\end{align*}
		where $K_2$ is sufficiently large, provided that either
		$\mathcal{H}_2$ or $\mathcal{H}_3$ holds. Indeed, it suffices to observe  that
		\begin{align*}
			\int_0^T b(t)[(1\ast b)(t)]^{-\tau(\beta,p)}\, \mathrm{d}t\lesssim [(1\ast b)]^{1-\tau(\beta,p)} \text{ ~for~ } \tau(\beta,p) \in [-1,1),
		\end{align*}
		 whence Young' inequality ($1+\frac{1}{q} = \frac{1}{1}+\frac{1}{q}$) yields this conclusion.
	\end{remark}
	
	\section{Averaging effects}
	Assume that \eqref{2.1} holds. Let $u^\omega \in L^2(\Omega, H^{s,p}(\mathbb R^3))$ be the random initial value corresponding to $u \in H^{s,p}(\mathbb R^3)$
	with $p \in (1,2)$ and $s\geq 0$.
	Consider the free evolution with data $u^\omega$, given by $w_u^\omega(t,x) = \mathcal C(t) u^\omega$. The subsequent averaging effects serve to estimate the probability of specific subsets in $\Omega$.
	\begin{proposition}\label{prop 5.1}
		Let $(a,b)\in (\mathcal{PC}^{\ast})$. Assume that $(\mathcal{H}_1)$ holds and $0\leq q \tau_1(\beta,p,s) < 1$ for some $q \in [1, \infty)$. For any fixed $T \in (0,1]$, we define the set
		\begin{align*}
			G_{\varsigma,q,p'} = \{\omega \in \Omega: \|w_u^\omega\|_{L^q([0,T];L^{p'}(\mathbb R^3))}\geq \varsigma \}.
		\end{align*}
		Then for each $\varsigma>0$, for each $u \in H^{s,p}(\mathbb R^3)$, there exist $A_1>0$ and $L_1>0$ such that
		\begin{align*}
			\mathbb P(G_{\varsigma, q, p'}) \leq L_1 e^{-A_1 \varsigma^2 / \|u\|^2_{H^{s,p}}}
		\end{align*}
		\begin{proof}
			It follows from the Minkowski's inequality, Lemma \ref{lemma 2.3} and \eqref{4.3} that
			\begin{align*}
				&\quad \|w_u^\omega\|_{L^{p_2}(\Omega;{L^{q}([0,T];L^{p'}(\mathbb R^3))})}\\
				&= \Bigg\|\sum_{j=0 }^{\infty} X_j(\omega) \mathcal F^{-1}( \varrho_j(\xi)C(t,\xi)\tilde{u}(\xi)) \Bigg\|_{L^{p_2}(\Omega;{L^{q}([0,T];L^{p'}(\mathbb R^3))})}\\
				&\leq \Bigg\| \Big\|\sum_{j=0}^{\infty} X_j(\omega)\mathcal F^{-1}( \varrho_j(\xi)C(t,\xi)\tilde{u}(\xi)) \Big\|_{L^{p_2}(\Omega)} \Bigg\|_{L^q([0,T];L^{p'}(\mathbb R^3))}\\
				&\lesssim \sqrt{p_2}\Bigg\| \Big(\sum_{j=0}^{\infty} \big| \mathcal F^{-1}( \varrho_j(\xi)C(t,\xi)\tilde{u}(\xi)) \big|^2 \Big)^{\frac{1}{2}} \Bigg\|_{L^q([0,T];L^{p'}(\mathbb R^3))}\\
				&\leq \sqrt{p_2} \Bigg\|\Big(\sum_{j=0}^{\infty} \big\| \mathcal F^{-1}( \varrho_j(\xi)C(t,\xi)\tilde{u}(\xi)) \big\|^2_{L^{p'}(\mathbb R^3) } \Big)^{\frac{1}{2}} \Bigg\|_{L^q([0,T])}\\
				&=\sqrt{p_2} \Bigg\|\Big(\sum_{j=0}^{\infty} \big\| \mathcal{C}(t) \mathcal F^{-1}(\varrho_j\tilde u) \big\|^2_{L^{p'}(\mathbb R^3)} \Big)^{\frac{1}{2}}\Bigg\|_{L^q([0,T])} \\
				&\lesssim  \Bigg(\int_0^T [m_0+(1\ast m_1)(t)]^{q\tau_1(\beta,p,s) } t^{-q\tau_1(\beta,p,s)}\Bigg)^{\frac{1}{q} } \sqrt{p_2} \Big(\sum_{j=0}^{\infty} \big\| \mathcal F^{-1}(\varrho_j\tilde u) \big\|^2_{H^{s,p}(\mathbb R^3)} \Big)^{\frac{1}{2}} \\
				&\leq \sqrt{p_2} M_1(T) \Bigg(\sum_{j=0}^{\infty} \big\|\mathcal F^{-1}(\varrho_j\tilde u)\big\|^2_{H^{s,p}(\mathbb R^3)} \Bigg)^{\frac{1}{2}} \\
				&\lesssim \sqrt{p_2} \|u\|_{H^{s,p}(\mathbb R^3)},
			\end{align*}
			provided that $p_2\geq \max\{q, p'\}$. Hence, invoking the Bienaym\'{e}-Tchebichev inequality, we obtain that there exists a constant $\alpha>0$ such that for every $p_2 \geq \max\{q, p'\}$, every $u \in H^{s,p}(\mathbb R^3)$, it holds that
			\begin{align}\label{5.1}
				\mathbb P(G_{\varsigma, q, p'}) \leq \varsigma^{-p_2} (\alpha \sqrt{p_2}\|u\|_{H^{s,p}})^{p_2}.
			\end{align}
			If $\varsigma$ satisfies $\frac{\varsigma}{\|u\|_{H^s}} \leq \max\{q, p'\}\alpha e$, then there exist some $L_1>1$ and $A_1\geq (\alpha e)^{-2}$ such that
			\begin{align*}
				L_1e^{-A_1 \varsigma^2 / \|u\|^2_{H^s} }\geq L_1e^{-(\max\{q, p'\}\alpha e A_1)^2}\geq 1\geq \mathbb P(G_{\varsigma, q, p'}).
			\end{align*}
			If $\varsigma$ satisfies $\frac{\varsigma}{\|u\|_{H^s}}>  \max\{q, p'\} \alpha e$, we put $p_2:= (\frac{\varsigma}{\alpha e \|u_0\|_{H^s} })^2 > (\max\{q, p'\})^2$. It follows from \eqref{5.1} that		
			\begin{align*}
				\mathbb P(E_{\varsigma, q, p'}) \leq \varsigma^{-(\frac{\varsigma}{\alpha e \|u\|_{H^s} })^2} (\frac{\varsigma}{e })^{(\frac{\varsigma}{\alpha e \|u\|_{H^s} })^2}=e^{ -(\alpha e)^{-2}\varsigma^2/ \|u\|^2_{H^s} }.
			\end{align*}
			This complete the proof.
		\end{proof}
	\end{proposition}
	
	\begin{proposition}\label{prop 5.2}
		Assume that $(a,b)\in (\mathcal{PC}^{\ast})$, $p \in (1,2)$ and $s\in \mathbb R_+$. Define the set
		\begin{align*}
			E_{\varsigma,p,s} = \{\omega \in \Omega: \|u^\omega\|_{H^{s,p}(\mathbb R^3)}\geq \varsigma \}.
		\end{align*}
		Then for each $\varsigma>0$, for each $u \in H^{s,p}(\mathbb R^3)$, there exist $A_2>0$ and $L_2>0$ such that
		\begin{align*}
			\mathbb P(E_{\varsigma,p,s}) \leq L_2 e^{-A_2 \varsigma^2 / \|u\|^2_{H^{s,p}}}
		\end{align*}
		\begin{proof}
			It follows from the Minkowski's inequality, Lemma \ref{lemma 2.3} that
			\begin{align*}
				\|u^\omega\|_{L^{p_2}(\Omega;{H^{s,p}(\mathbb R^3)})}
				&= \Bigg\|\sum_{j=0 }^{\infty} X_j(\omega) \mathcal F^{-1}( \varrho_j \tilde{u}) \Bigg\|_{L^{p_2}(\Omega;H^{s,p}(\mathbb R^3))}\\
				&= \Bigg\|\sum_{j=0 }^{\infty} X_j(\omega) \mathcal F^{-1}( \varrho_j(\xi) (1+|\xi|^2)^{\frac{s}{2} } \tilde{u}(\xi)) \Bigg\|_{L^{p_2}(\Omega;L^{p}(\mathbb R^3))}\\
				&\leq \Bigg\| \Big\|\sum_{j=0}^{\infty} X_j(\omega)\mathcal F^{-1}( \varrho_j(\xi) (1+|\xi|^2)^{\frac{s}{2} } \tilde{u}(\xi)) \Big\|_{L^{p_2}(\Omega)} \Bigg\|_{L^{p}(\mathbb R^3)}\\
				&\lesssim \sqrt{p_2}\Bigg\| \Big(\sum_{j=0}^{\infty} \big| \mathcal F^{-1}( \varrho_j(\xi) (1+|\xi|^2)^{\frac{s}{2} }\tilde{u}(\xi)) \big|^2 \Big)^{\frac{1}{2}} \Bigg\|_{L^{p}(\mathbb R^3)}\\
				&\lesssim \sqrt{p_2} \|u\|_{H^{s,p}(\mathbb R^3)},
			\end{align*}
			provided that $p_2\geq 2$. Hence, invoking the Bienaym\'{e}-Tchebichev inequality, we obtain that there exists a constant $\alpha>0$ such that for every $p_2 \geq 2$ and every $u \in H^s(\mathbb R^3)$,  the following estimate holds:
			\begin{align*}
				\mathbb P(E_{\varsigma,p,s}) \leq \varsigma^{-p_2} (\alpha \sqrt{p_2}\|u\|_{H^s})^{p_2}.
			\end{align*}
			The following proof is similar to that in Proposition \ref{prop 5.1}, so we omit it.
		\end{proof}
	\end{proposition}
	
	\section{Existence and uniqueness for the random initial data}
	In this section, we establish the almost sure existence of local solutions for randomized initial data. For the purpose of solving the problem
	\begin{align}\label{6.1}
		\begin{cases}
			\partial_t^2 \left(\mathrm{d}a \ast \mathrm{d}a \ast (w-u^\omega)\right)(t) +\gamma \partial_t \left(\mathrm{d}a \ast (w-u^\omega) \right)(t) +(- \Delta)^{\frac{\beta}{2}}w =- w|w|^{\kappa-1}, \\
			w_{t=0}=u^\omega(x), \ \partial_t w_{t=0}=0,
		\end{cases}
	\end{align}
	let $w = w^\omega_u + v$, we conclude that $v$ solves
	\begin{align*}
		\begin{cases}
			\partial_t^2 \left(\mathrm{d}a \ast \mathrm{d}a \ast v \right)(t) +\gamma \partial_t \left(\mathrm{d}a \ast v \right)(t) +(- \Delta)^{\frac{\beta}{2}}v =- (v+w^\omega_u)|v+w^\omega_u|^{\kappa-1}, \\
			v_{t=0}=0, \ \partial_t v_{t=0}=0.
		\end{cases}
	\end{align*}
	In other words, $v(t) = -\int_0^t \mathcal S(t-\zeta) (v+w^\omega_{u})|v+w^\omega_{u} |^{\kappa-1}\, \mathrm{d}\zeta$. Define the map
	\begin{align*}
		W^\omega_u: v \mapsto -\int_0^t \mathcal S(t-\zeta) (v+w^\omega_{u})|v+w^\omega_{u} |^\kappa\, \mathrm{d}\zeta.
	\end{align*}
	We present the local existence below. Our analysis requires the following parameter constraints for $\beta, p$ and $s$.
	\begin{enumerate}
		\item[$(\mathcal{H}_4)$] $\kappa \in (2,3]$ and $s \in [\frac{2(\kappa-1)}{\kappa+1}, \frac{3(\kappa-1)}{\kappa+1})$ if $\beta =2$, and $\kappa \in (\frac{6+\beta}{6-\beta},\frac{12-\beta}{12-5\beta} ]$ and $s \in (0, 2-\frac{\beta}{2})\cap [\frac{(12-3\beta)(\kappa-1)}{4(\kappa+1)}, \frac{3(\kappa-1)}{\kappa+1})$ if $\beta \in (1,2)$.
	\end{enumerate}
	
	\begin{enumerate}
		\item[$(\mathcal{H}_5)$] $\kappa \in (1, 2]$ and $s \in [\frac{2(\kappa-1)}{\kappa+1}, \frac{3(\kappa-1)}{\kappa+1}]$ if $\beta =2$, and $\kappa \in (1, \frac{6+\beta}{6-\beta}]$ and $s \in (0, 2-\frac{\beta}{2})\cap [\frac{(12-3\beta)(\kappa-1)}{4(\kappa+1)}, \frac{3(\kappa-1)}{\kappa+1}]$ if $\beta \in (1,2)$.
	\end{enumerate}
	
	\begin{remark}\label{rem 6.1}
		Let $\beta \in (1,2)$. A simple calculation derives that
		\begin{align*}
			\begin{cases}
				\sup\limits_{\kappa \in (\frac{6+\beta}{6-\beta},\frac{12-\beta}{12-5\beta} ] } \frac{(12-3\beta)(\kappa-1)}{4(\kappa+1)}= \frac{\beta}{2}<2-\frac{\beta}{2},   \\
				\sup\limits_{\kappa \in (1, \frac{6+\beta}{6-\beta}] } \frac{(12-3\beta)(\kappa-1)}{4(\kappa+1)} = \frac{\beta}{4}(2-\frac{\beta}{2})<2-\frac{\beta}{2},
			\end{cases}
		\end{align*}
		which implies that the set $(0, 2-\frac{\beta}{2})\cap [\frac{(12-3\beta)(\kappa-1)}{4(\kappa+1)}, \frac{3(\kappa-1)}{\kappa+1})$ is not empty.
	\end{remark}
	
	\begin{theorem}\label{the 6.1}
		Assume that $(a,b)\in (\mathcal{PC}^{\ast})$, $b$ is  nonincreasing, $u^\omega \in L^2(\Omega, H^{s,1+\frac{1}{\kappa}}(\mathbb R^3))$ is the random initial value corresponding to $u \in H^{s,1+\frac{1}{\kappa}}(\mathbb R^3)$ and \eqref{2.1} holds. Let $q_0 \in (1,\kappa]$ and take $q = \frac{q_0(\kappa-1)}{q_0-1}$.
		If one of the following holds:
		\begin{itemize}
			\item[(i)] $a_0=0$, $(\mathcal{H}_4)$ holds, $b \in L_{loc}^{p_0}(\mathbb R_+)$ for some $p_0>\frac{q_0(1-\tau(\beta,p))}{1-q_0\tau(\beta,p)}$ with $0<q_0\tau(\beta,p)<1$, and $0<q\tau_1(\beta,p,s)<1$.
			\item[(ii)] $a_0=0$, $(\mathcal{H}_5)$ holds, $b \in L_{loc}^{q_0}(\mathbb R_+)$, and $0<q\tau_1(\beta,p,s)<1$.
			\item[(iii)] $a_0>0$, $(\mathcal{H}_4)$ holds,  $0<q_0\tau(\beta,p)<1$ and $0<q\tau_1(\beta,p,s)<1$.
			\item[(iv)] $a_0>0$, $(\mathcal{H}_5)$ holds and $0<q\tau_1(\beta,p,s)<1$.
		\end{itemize}
		Then for almost all $\omega \in \Omega$, there exist $T_{\omega} \in (0,1]$ and a unique local solution $w\in w_u^\omega + L^{q}([0,T_\omega];L^{1+\kappa}(\mathbb R^3))$
		to problem \eqref{6.1}. More precisely, for every $T \in (0,1]$, 	
		there exist $B_1>0$ and an event $\Omega_{1,T}$ satisfying
		\begin{align*}
			\mathbb P(\Omega_{1,T}) \geq 1- L_1 e^{-B_1 M(T)^{-\frac{2}{\kappa-1}} \big/\big.  { \big\|u \big\|^2_{H^{s,1+\frac{1}{\kappa} }}} },
		\end{align*}
		such that for every $\omega \in \Omega_{1,T}$, there exists a unique local mild solution $w \in w_u^\omega + L^{q}([0,T];L^{1+\kappa}(\mathbb R^3))$
		to problem \eqref{6.1}.
		\begin{proof}
			Let $C_{\kappa}:=\kappa 2^{\kappa-1}$ and $p=1+\frac{1}{\kappa}$.
			A simple calculation shows that $(\mathcal{H}_2)$ holds under assumption $(\mathcal{H}_4)$,  $(\mathcal{H}_3)$ holds under assumption $(\mathcal{H}_5)$, while $(\mathcal{H}_1)$ remains valid under either assumption. Consequently,
			Lemma \ref{lemma 4.2}(1)-(4) hold under Theorem \ref{the 6.1}(i)-(iv), respectively. Let
			$r=\frac{q_0(\kappa-1)}{(q_0-1)\kappa}$,
			then $q,r$ satisfy $1+\frac{1}{q} = \frac{1}{q_0} +\frac{1}{r}$ and $q = \kappa r$. Given $T\in (0,1]$ and $\omega \in G^c_{\varsigma,q,p'}$, Lemma \ref{lemma 4.2} yields that
			\begin{align*}
				\begin{aligned}
					\|W^\omega_u(v) \|_{L^q([0,T]; L^{p'})} &\leq M(T) \|(v+w^\omega)|v+w^\omega_u|^{\kappa-1}\|_{L^r([0,T]; L^{p})}\\
					&\leq M(T) \|v+w^\omega_u \|^\kappa_{L^q([0,T]; L^{p'})}\\
					&\leq 2^{\kappa-1}M(T) \left(\|v\|^\kappa_{L^q([0,T]; L^{p'})} + \|w^{\omega}_u\|^\kappa_{L^q([0,T]; L^{p'})} \right)\\
					&\leq C_{\kappa} M(T) \left(\|v\|^\kappa_{L^q([0,T]; L^{p'})} +  \varsigma^\kappa \right).
				\end{aligned}
			\end{align*}
			and
			\begin{align*}
				\begin{aligned}
					&\quad \left\|W^\omega_u (v_1) - W^\omega_u (v_2) \right\|_{L^q([0,T];L^{p'})}\\
					&\leq M(T) \left\|(v_1+ w^\omega_u)|v_1+w^\omega_u|^{\kappa-1} - (v_2+w^\omega_u)|v_2+w^\omega_u|^{\kappa-1} \right\|_{L^r([0,T];L^{p})}\\
					&\leq \kappa M(T) \left\|(v_1- v_2)\left(|v_1+w^\omega_u|^{\kappa-1}+|v_2+w^\omega_u|^{\kappa-1} \right)\right\|_{L^r([0,T];L^{p})}\\
					&\leq \kappa 2^{\kappa-1}M(T) \|v_1-v_2 \|_{L^q([0,T];L^{p'})}\\
					&\quad \times \left(2\|w^{\omega}_u\|^{\kappa-1}_{L^q([0,T]; L^{p'})}  +\|v_1\|^{\kappa-1}_{L^q([0,T];L^{p'})}+ \|v_2\|^{\kappa-1}_{L^q([0,T];L^{p'})} \right)\\
					&\leq C_{\kappa} M(T) \|v_1-v_2 \|_{L^q([0,T];L^{p'})}\left(2\varsigma^{\kappa-1}  +\|v_1\|^{\kappa-1}_{L^q([0,T];L^{p'})}+ \|v_2\|^{\kappa-1}_{L^q([0,T];L^{p'})} \right).
				\end{aligned}
			\end{align*}
			Take $\varsigma$ satisfying $C_{\kappa} M(T)\varsigma^{\kappa-1}< \frac{1}{6}$, then the map $W_u^\omega$ is a contraction on the ball of radius $\varsigma$ of $L^q([0,T]; L^{p'})$. Take $\varsigma(T) = [7C_{\kappa} M(T)]^{-\frac{1}{\kappa-1}}$, and set
			\begin{align*}
				\Omega_{1,T} = G^c_{\varsigma(T), q,p'}, \ \Phi_1 = \bigcup_{n \in \mathbb N_+} \Omega_{1,\frac{1}{n}}.
			\end{align*}
			Since $(\mathcal{H}_1)$ holds and $0<q\tau_1(\beta,p,s)<1$, by Proposition \ref{prop 5.1} we obtain
			\begin{align*}
				\mathbb P(\Omega_{1,T}) \geq 1- L_1 e^{-B_1 M(T)^{-{\frac{2}{\kappa-1}}} / \|u\|^2_{H^{s,p}}}, \ \mathbb P(\Phi_1) = 1 ,
			\end{align*}
			where $B_1 = A_1(7C_{\kappa})^{-\frac{2}{\kappa-1}}$.
			The proof is completed.
		\end{proof}
	\end{theorem}
	
	We now present two representative examples to illustrate the above results.
	
	\begin{example}
		Damped cubic wave equation:
		\begin{align}\label{6.2}
			\begin{cases}
				\partial^2_t w + \partial_t w -\Delta w =- w|w|^{2}, \quad t>0, \ x \in \mathbb R^3, \\
				w_{t=0}=u^\omega(x), \ \partial_t w_{t=0}=0,  \quad  x \in \mathbb R^3.
			\end{cases}
		\end{align}
		The kernel $a$ is the Heaviside function, which means $a_0 =1$, $a_1 \equiv 0$. Take $b$ $\equiv 1$, then $\mathrm{d}a \ast b \equiv 1$. Thus, $(\delta_0, 1) \in (\mathcal{PC}^{\ast})$. If $u \in H^{\frac{11}{8},\frac{4}{3}}$, invoking Theorem \ref{the 6.1}, for every $T \in (0,1]$, there exists an event $\Omega_{1,T}$ satisfying
		\begin{align*}
			\mathbb P(\Omega_{1,T}) \geq 1- L_1 e^{-\frac{B_1}{K_2} T^{-\frac{1}{10}} \big / \big. \big\|u \big \|^2_{H^{\frac{11}{8},\frac{4}{3}  }}  },
		\end{align*}
		such that for almost every $\omega \in \Omega_{1,T}$, there exists a unique local mild solution $w \in w_u^\omega + L^{5}([0,T];L^{4}(\mathbb R^3))$
		to problem \eqref{6.2}, where
		\begin{align*}
			\Omega_{1,T} = \{\omega \in \Omega: \|w_u^\omega\|_{L^5([0,T];L^{4}(\mathbb R^3))}\leq  (7C_3)^{-\frac{1}{2}}T^{-\frac{1}{20} } \}.
		\end{align*}
		Here, setting $q_0 =\frac{5}{3}$ and $s=\frac{11}{8}$, straightforward computation yields
		\begin{align*}
			p=\frac{4}{3}, \ q = 5, \ q_0\tau(2,p) = \frac{5}{6}, \ q\tau_1(2,p,s) = \frac{5}{8}, \ M(T) =K_2 T^{\frac{1}{10} }.
		\end{align*}
		Thus, the condition Theorem \ref{the 6.1}(iii) holds.
	\end{example}
	
	\begin{example}
		The space-time fractional telegraph equation:
		\begin{align}\label{6.3}
			\begin{cases}
				\partial_t^{\frac{11}{6}} w(t) +\gamma \partial_t^{\frac{11}{12}} w(t)  +(- \Delta)^{\frac{2}{3}}w =- w|w|, \quad t>0, \ x \in \mathbb R^3, \\
				w_{t=0}=u^\omega(x), \ \partial_t w_{t=0}=0,  \quad  x \in \mathbb R^3,
			\end{cases}
		\end{align}
		where $\gamma>0$ and $\partial_t^\alpha$ represents the Caputo fractional derivative. The kernel $a = 1 \ast g_{1-\alpha}$, which means $a_0 =0$, $a_1=g_{1-\alpha}$ and $b = g_\alpha$. If $u \in H^{\frac{5}{6},\frac{3}{2}}$, invoking Theorem \ref{the 6.1}, for every $T \in (0,1]$, there exists an event $\Omega_{1,T}$ satisfying
		\begin{align*}
			\mathbb P(\Omega_{1,T}) \geq 1- L_1 e^{-\frac{B_1}{K_2} T^{-\frac{1}{10}} \big/\big. \big\|u \big \|^2_{H^{\frac{5}{6},\frac{3}{2}  }} }    ,
		\end{align*}
		such that for every $\omega \in \Omega_{1,T}$, there exists a unique local mild solution $w \in w_u^\omega + L^{\frac{5}{2}}([0,T];L^{3}(\mathbb R^3))$
		to problem \eqref{6.3}, where
		\begin{align*}
			\Omega_{1,T} = \{\omega \in \Omega: \|w_u^\omega\|_{L^{\frac{5}{2} }([0,T];L^{3}(\mathbb R^3))}\leq  (7C_2)^{-\frac{1}{2}}T^{-\frac{1}{20} } \}.
		\end{align*}
		Here, setting $q_0 =\frac{5}{3}$ and $s=\frac{5}{6}$, straightforward computation yields
		\begin{align*}
			p=\frac{3}{2}, \ q =\frac{5}{2}, \ q_0\tau(\frac{4}{3},p) = \frac{5}{6}, \ q\tau_1(\frac{4}{3},p,s) = \frac{5}{8}, \ M(T) = K_2 T^{\frac{1}{10} }.
		\end{align*}
		Thus, the condition Theorem \ref{the 6.1}(i) holds.
	\end{example}
	
	For the case where $s\geq \frac{3(\kappa-1)}{\kappa+1}$,
	we can prove that problem \eqref{6.1} admits one solution with enhanced temporal regularity.
	\begin{theorem}\label{the 6.2}
		Assume that $(a,b)\in (\mathcal{PC}^{\ast})$ and $b$ is nonincreasing, $u^\omega \in L^2(\Omega, H^{s,1+\frac{1}{\kappa}}(\mathbb R^3))$ is the random initial value corresponding to $u \in H^{s,1+\frac{1}{\kappa}}(\mathbb R^3)$ with $s \in [ \frac{3(\kappa-1)}{\kappa+1}, \infty)$ and \eqref{2.1} holds. Then for almost all $\omega \in \Omega$, there exist $T_{\omega} \in (0,1]$ and a unique solution $w\in w_u^\omega + L^{\infty}([0,T_\omega];L^{1+\kappa}(\mathbb R^3))$ to problem \eqref{6.1} provided
		\begin{align}\label{6.4}
			\begin{cases}
				\kappa \in (1,3] \text{ ~if~ } \beta =2,\\
				\kappa \in  (1, \frac{12-\beta}{12-5\beta}] \cap (1, \frac{10-\beta}{2+\beta}) \text{ ~if~ }\beta \in (1,2).
			\end{cases}
		\end{align}
		More precisely, for every $T \in (0,1]$, there exist $B_0>0$ and  an event $\Omega_{0,T}$ satisfying
		\begin{align*}
			\mathbb P(\Omega_{0,T}) \geq 1- L_2 e^{-B_0  \|b\|_{L^1([0,\frac{T}{4}])}^{-{\frac{2(1-\tau(\beta,p) )}{\kappa-1}}} \big/\big. \big\|u \big\|^2_{H^{s,1+\frac{1}{\kappa} }}},
		\end{align*}
		such that for every $\omega \in \Omega_{0,T}$, there exists a unique local mild solution $w \in w_u^\omega + L^{\infty}([0,T];L^{1+\kappa}(\mathbb R^3))$ to problem \eqref{6.1}.
		\begin{proof}
			Let $s_0= \frac{3(\kappa-1)}{\kappa+1}$ and $p=1+\frac{1}{\kappa}$. A straightforward calculation shows that depending on the values of $\kappa$ and $\beta$ under condition \eqref{6.4}, either $(\mathcal{H}_2)$ or $(\mathcal{H}_3)$ will be satisfied.
			From the relation $s_0 = \frac{3(\kappa-1)}{\kappa+1}$, it immediately follows that $\tau_1(\beta,p,s_0) = 0$ and $p =\frac{6}{3+s_0}$. Note that \eqref{6.4} yields
			\begin{align*}
				\begin{cases}
					s_0 \in (0,\frac{3}{2}]  \text{ ~if~ } \beta =2,\\
					s_0 \in  (0, \frac{2\beta}{4-\beta}] \cap (0, 2-\frac{\beta}{2})  \text{ ~if~ }\beta \in (1,2),
				\end{cases}
			\end{align*}
			which shows that $\beta$, $p$, and $s_0$ (which replaces $s$) together satisfy hypothesis $(\mathcal{H}_1)$. Note that $H^{s,p} \hookrightarrow H^{s_0,p}$ indicates that there exist a constant $I_{s,s_0}$ such that $\|u^\omega\|_{H^{s_0,p}} \leq I_{s,s_0} \|u^\omega\|_{H^{s,p}}$.
			Take $C_{\kappa,1} =2^{\kappa-1}\max\{\kappa , (I_{s,s_0}K_1)^\kappa, (I_{s,s_0}K_1)^{\kappa-1} \}$. Given $T \in (0,1]$ and $\omega \in E^c_{\varsigma,p,s}$, one can derive from Remark \ref{rem 4.2} and Lemma \ref{lemma 4.1} that
			\begin{align*}
				\|W^\omega_u(v) \|_{L^\infty([0,T]; L^{p'})}
				&\leq 2^{\kappa-1}K_2 \|b\|^{1-\tau(\beta,p)}_{L^1([0,\frac{T}{4}])} \left(\|v\|^\kappa_{L^\infty([0,T]; L^{p'})} +  \|\mathcal C(t)u^\omega\|^\kappa_{L^\infty([0,T]; L^{p'})} \right)\\
				&\leq 2^{\kappa-1}K_2 \|b\|^{1-\tau(\beta,p)}_{L^1([0,\frac{T}{4}])} \left(\|v\|^\kappa_{L^\infty([0,T]; L^{p'})} +  K_1^\kappa\|u^\omega\|^\kappa_{H^{s_0,p}} \right)\\
				&\leq 2^{\kappa-1}K_2 \|b\|^{1-\tau(\beta,p)}_{L^1([0,\frac{T}{4}])} \left(\|v\|^\kappa_{L^\infty([0,T]; L^{p'})} +  (I_{s,s_0}K_1)^\kappa\|u^\omega\|^\kappa_{H^{s,p}} \right)\\
				&\leq C_{\kappa,1} K_2 \|b\|^{1-\tau(\beta,p)}_{L^1([0,\frac{T}{4}])} \left(\|v\|^\kappa_{L^\infty([0,T]; L^{p'})} +  \varsigma^\kappa \right)
			\end{align*}
			and
			\begin{align*}
				&\quad \left\|W^\omega_u (v_1) - W^\omega_u (v_2) \right\|_{L^\infty([0,T];L^{p'})}\\
				&\leq C_{\kappa,1} K_2 \|b\|^{1-\tau(\beta,p)}_{L^1([0,\frac{T}{4}])} \|v_1-v_2 \|_{L^\infty([0,T];L^{p'})}\left(2\varsigma^{\kappa-1}  +\|v_1\|^{\kappa-1}_{L^\infty([0,T];L^{p'})}+ \|v_2\|^{\kappa-1}_{L^\infty([0,T];L^{p'})} \right).
			\end{align*}
			Take $\varsigma(T) = [7C_{\kappa,1} K_2 \|b\|^{1-\tau(\beta,p)}_{L^1[0,\frac{T}{4}]}]^{-\frac{1}{\kappa-1}}$. Then the map $W_u^\omega$ is a contraction on the ball of radius $\varsigma$ of $L^\infty([0,T]; L^{p'})$. Set $\Omega_{0,T} = E^c_{\varsigma(T),p,s}$ and $\Phi_0 = \bigcup_{n \in \mathbb N_+} \Omega_{0,\frac{1}{n}}$. Finally, by Proposition \ref{prop 5.2} we obtain
			\begin{align*}
				\mathbb P(\Omega_{0,T}) \geq 1- L_2 e^{-B_0  \|b\|_{L^1([0,\frac{T}{4}])}^{-{\frac{2(1-\tau(\beta,p) )}{\kappa-1}}} / \|u\|^2_{H^{s,p}}}, \ \mathbb P(\Phi_0) = 1 ,
			\end{align*}
			where $B_0 = A_2(7C_{\kappa,1}K_2)^{-\frac{2}{\kappa-1}}$. The proof is completed.
		\end{proof}
	\end{theorem}
	
	\begin{remark}\label{rem 6.2}
		Theorem \ref{the 6.1} and Theorems \ref{the 6.2} reveal a critical index phenomenon characterized by the critical regularity index
		$s_{crit} = \frac{3(\kappa-1)}{\kappa+1}$. For randomized initial data, the solution admits
		$L_t^q$-regularity ($q<\infty$) with probability one. Conversely, for any $s\geq s_{crit}$, the solution attains $L_t^\infty$-regularity almost surely.
	\end{remark}
	
	\begin{appendices}
	\section{Auxiliary Lemma concerning relaxation functions}\label{app a}
	The following lemmas constitute an improvement and generalization of existing results.
	Note that $k_1$ is nonnegative nonincreasing in $(0,\infty)$ when $l$ is completely positive. Let $k_\infty:=\lim\limits_{t \to \infty} k_1(t)$, then there exists a nonnegative nonincreasing function $k_2\in L^1_{loc}(\mathbb R_+)$ such that $k_1(t)=k_\infty +k_2(t)$ for $t>0$. Motivated by the works in \cite{P. Clement 79, J.C. Pozo}, we generalize and obtain the following lemma.
	\begin{lemma}\label{lemma A.1}
		Assume that $l\in L^1_{loc}(\mathbb R_+)$ is completely positive.
		Let $\gamma \in \mathbb C$. $s_\gamma(t), r_\gamma(t)$ are defined by \eqref{2.4} and \eqref{2.5}, respectively. Then we have the following results.
		\begin{enumerate}
			\item[(i)]
			For each $t\geq 0$, there exists a positive finite measure $-\phi(t, \mathrm{d}\theta)$ such that
			\begin{align*}
				-\widehat{\phi}(\lambda, \mathrm{d}\theta) = (\lambda \widehat{l}(\lambda))^{-1} e^{-\theta \widehat {l}(\lambda)^{-1}}\mathrm{d}\theta.
			\end{align*}
			Furthermore, for $\gamma \in \{z \in \mathbb {C}: Re(z)\geq -k_\infty\} $, $s_\gamma(t)$ admits the representation
			\begin{align*}
				s_\gamma(t) = -\int_0^\infty e^{-\gamma \theta} \phi(t,\mathrm{d}\theta), \quad t\geq 0.
			\end{align*}
			Particularly, setting $\gamma =0$ yields $-\int_0^\infty \phi(t,\mathrm{d}\theta)=1, \ t\geq 0$.
			\item[(ii)] For each $t>0$, there exists a positive finite measure $\eta(t, \mathrm{d}\theta)$ such that
			\begin{align*}
				\widehat{\eta}(\lambda, \mathrm{d}\theta) = e^{-\theta \widehat{l}(\lambda)^{-1}}\mathrm{d}\theta.
			\end{align*}
			Furthermore, for $\gamma \in \{z \in \mathbb {C}: Re(z)\geq -k_\infty\} $, $r_\gamma(t)$ admits the representation
			\begin{align*}
				r_\gamma(t) = \int_0^\infty e^{-\gamma \theta} \eta(t,\mathrm{d}\theta), \quad \ t>0.
			\end{align*}
			Particularly, setting $\gamma =0$ yields $\int_0^\infty \eta(t,\mathrm{d}\theta)=l(t), \ t>0$.
		\end{enumerate}
		\begin{proof}
			Let us first prove (i). It is shown in \cite[Proposition 4.5]{J. Pruss} that $(k,l) \in (\mathcal \mathcal{PC^{\ast}})$ is equivalent to
			$[{\widehat{l}(\lambda)}]^{-1} \in (\mathcal{BF})$. Further, $\lambda \mapsto e^{-{\theta}{\widehat{l}(\lambda)}^{-1}}\in (\mathcal{CM})$ for each $\theta \geq 0$, see \cite[Proposition 4.2]{J. Pruss}. Applying Bernstein's Theorem (\cite[Chapter 4]{J. Pruss}), there exists a unique function $\phi(\cdot, \theta) \in BV_{loc}(\mathbb R_+)$
			that is nondecreasing with $\phi(0,\theta)=0$ such that
			\begin{align}\label{A.1}
				\widehat{\mathrm{d}_t \phi}(\lambda, \theta) = e^{-\theta \widehat{l}(\lambda)^{-1}}, \quad \theta \geq 0.
			\end{align}
			It follows immediately that
			$\widehat{\phi}(\lambda, \theta) = {\lambda}^{-1}e^{-\theta \widehat{l}(\lambda)^{-1}}$. Differentiating both sides with respect to $\theta$ leads to
			\begin{align}\label{A.2}
				-\widehat{\phi}(\lambda, \mathrm{d}\theta) = (\lambda \widehat{l}(\lambda))^{-1}e^{-\theta \widehat{l}(\lambda)^{-1}}\mathrm{d}\theta,
			\end{align}
			According to \eqref{2.6}, we obtain $\widehat {l}(\lambda)^{-1} = k_0\lambda+\lambda\widehat{k_1}(\lambda)=k_\infty+k_0\lambda+\lambda\widehat{k_2}(\lambda)$ for $\lambda>0$. It is stated in \cite[Proposition 4.3]{J. Pruss} that $\inf\limits_{\lambda>0} \widehat{l}(\lambda)^{-1}=k_\infty$. Thus, when $\gamma \in \{z \in \mathbb {C}: Re(z)\geq -k_\infty\} $, Integrating \eqref{A.2} yields
			\begin{align}\label{A.3}
				-\int_0^\infty e^{-\gamma \theta} \widehat{\phi}(\lambda, \mathrm{d}\theta) = \frac{1}{\lambda}\frac{1}{1+\gamma \widehat{l}(\lambda)}=\widehat{s_\gamma}(\lambda), \quad Re(\lambda)>0.
			\end{align}
			Due to the uniqueness of Laplace transforms, statement (i) holds.
			
			Now turn to (ii). It follows from \eqref{A.1} that
			\begin{align*}
				\widehat{l}(\lambda)e^{-\theta \widehat{l}(\lambda)^{-1}}= \widehat{l}(\lambda) \widehat{\mathrm{d}_t \phi}(\lambda,\theta) = \widehat{l \ast \mathrm{d}_t \phi}(\lambda, \theta).
			\end{align*}
			Differentiating both sides with respect to $\theta$ results in
			\begin{align*}
				e^{-\theta \widehat{l}(\lambda)^{-1}}\mathrm{d}\theta= - \widehat{l \ast \mathrm{d}_t \phi}(\lambda, \mathrm{d}\theta) := \widehat{\eta}(\lambda, \mathrm{d}\theta),
			\end{align*}
			Similar to the derivation of \eqref{A.3}, for $\gamma \in \{z \in \mathbb {C}: Re(z)\geq -k_\infty\} $, it holds that
			\begin{align*}
				\int_0^\infty e^{-\gamma \theta} \widehat{\eta}(\lambda, \mathrm{d}\theta) = \frac{\widehat{l}(\lambda)}{1+\gamma \widehat{l}(\lambda)}= \widehat{r_\gamma}(\lambda), \quad Re(\lambda)>0,
			\end{align*}
			which means (ii) holds.
		\end{proof}
	\end{lemma}
	The lemma below serves to address the $L_tL_x$ estimates to problem \eqref{1.4}. The proof is a direct adaptation of that for \cite[Lemma 2, Lemma 3]{J. C. Pozo 24}. The key modification lies in that we relax the condition $(\mathcal {PC})$ to the complete positive condition.
	\begin{lemma}\label{lemma A.2}
		Assume that $l\in L^1_{loc}(\mathbb R_+)$ is completely positive. Then, for each $\delta \in (0,1)$, there exist two nonnegative functions $c_\delta, d_\delta \in L_{loc}^1(\mathbb R_+)$, whose Laplace transforms are given by
		\begin{align*}
			\widehat{c}_\delta(\lambda) = \Gamma(1-\delta)\lambda^{-1}\widehat{l}(\lambda)^{-\delta}, \
			\widehat{d}_\delta(\lambda) = \Gamma(1-\delta)\widehat{l}(\lambda)^{1-\delta}, \quad Re(\lambda)>0,
		\end{align*}
		respectively, and $c_\delta$, $d_\delta$ admit the representation
		\begin{align*}
			&c_\delta(t)=\frac{1}{\Gamma(\delta)} \int_0^\infty \, \gamma^{\delta-1}s_\gamma(t) \mathrm{d}\gamma = -\int_0^\infty \theta^{-\delta} \phi(t, \mathrm{d} \theta), \\
			&d_\delta(t)=\frac{1}{\Gamma(\delta)} \int_0^\infty \, \gamma^{\delta-1}r_\gamma(t) \mathrm{d}\gamma = \int_0^\infty \theta^{-\delta} \eta(t, \mathrm{d} \theta),  \quad t>0,
		\end{align*}
		respectively. Moreover, it holds that
		\begin{align*}
			c_\delta(t) \leq \Gamma(1-\delta) [(1\ast l)(t)]^{-\delta}, \quad t>0,
		\end{align*}
		where $\Gamma(\cdot)$ stands for Gamma function. If $l$ is also nonincreasing, then
		\begin{align*}
			d_\delta(t) \leq \Gamma(1-\delta) l(t)[(1\ast l)(t)]^{-\delta}, \quad t>0.
		\end{align*}
		\begin{proof}
			We begin by establishing the existence of $c_\delta$. From Proposition \ref{prop 2.1}(ii) we obtain
			\begin{align*}
				\frac{1}{\Gamma(\delta)}\int_0^\infty \gamma^{\delta-1} s_\gamma(t) \, \mathrm{d}\gamma  \leq \frac{1}{\Gamma(\delta)}\int_0^\infty \frac{\gamma^{\delta-1}}{1+\gamma(1\ast l)(t)} \, \mathrm{d}\gamma = \Gamma(1-\delta) [(1\ast l)(t)]^{-\delta},\quad t>0.
			\end{align*}
			Here, we employ the fact
			\begin{align*}
				\frac{1}{\Gamma(\delta)}\int_0^\infty \gamma^{\delta-1}(1+\gamma)^{-1}\, \mathrm{d}\gamma = \frac{B(\delta, 1-\delta)}{\Gamma(\delta)} = \Gamma(1-\delta),
			\end{align*}
			where $B(\cdot,\cdot)$ denote the Beta function. Further, Lemma \ref{lemma A.1}(i) along with Fubini's theorem shows that
			\begin{align*}
				\int_0^\infty \gamma^{\delta-1} s_\gamma(t) \, \mathrm{d}\gamma  &= -\int_0^\infty \gamma^{\delta-1} \int_0^\infty e^{-\gamma \theta} \phi(t,\mathrm{d}\theta)\mathrm{d}\gamma =-\int_0^\infty \phi(t,\mathrm{d}\theta) \int_0^\infty \gamma^{\delta-1} e^{-\gamma \theta}\, \mathrm{d}\gamma\\
				&=\Gamma(\delta) \int_0^\infty \theta^{-\delta} \phi(t,\mathrm{d}\theta).
			\end{align*}
			Define the function $c_\delta(t): = \Gamma(\delta)^{-1}\int_0^\infty \gamma^{\delta-1}s_\gamma(t)\, \mathrm{d}\gamma$. Obviously, $c_\delta$ is nonnegative since $s_\gamma$ is nonnegative. Note that \eqref{2.6} indicates
			\begin{align*}
				[(1 \ast l)(t)]^{-1} \leq \frac{k_0+(1 \ast k_1)(t)}{t}, \quad t>0.
			\end{align*}
			Immediately, $c_\delta \lesssim t^{-\delta}$ for each fixed $T>0$, which implies $c_\delta \in L_{loc}^1(\mathbb R_+)$.
			
			Similarly, invoking Lemma \ref{lemma A.1}(ii) together with Fubini's theorem, we conclude that
			\begin{align*}
				\int_0^\infty \gamma^{\delta-1} r_\gamma(t) \, \mathrm{d}\gamma  &= \int_0^\infty \gamma^{\delta-1} \int_0^\infty e^{-\gamma \theta} \eta(t,\mathrm{d}\theta)\mathrm{d}\gamma =\int_0^\infty \eta(t,\mathrm{d}\theta) \int_0^\infty \gamma^{\delta-1} e^{-\gamma \theta}\, \mathrm{d}\gamma\\
				&=\Gamma(\delta) \int_0^\infty \theta^{-\delta} \eta(t,\mathrm{d}\theta).
			\end{align*}
			Define the nonnegative function $d_\delta(t): =\Gamma(\delta)^{-1} \int_0^\infty \gamma^{\delta-1}r_\gamma(t)\, \mathrm{d}\gamma$ for $t>0$. For every fixed $T>0$, by Proposition \ref{prop 2.1}(ii) we get
			\begin{align*}
				\gamma^{\delta-1} \int_0^T r_\gamma(t)\, \mathrm{d}t =  \gamma^{\delta-2}(1-s_\gamma(T))
				\leq  \gamma^{\delta-1} \frac{[k_1(T)]^{-1}}{1+[k_1(T)]^{-1}\gamma}, \quad \gamma>0,
			\end{align*}
			where the right of the above inequality belong to $L_\gamma^1([0,T])$, then Fubini's  theorem yields $d_{\delta} \in L^1_{loc}(\mathbb R_+)$. If $l$ is also nonincreasing, it follows from Proposition \ref{prop 2.1}(iii) that
			\begin{align*}
				\int_0^\infty \gamma^{\delta-1} r_\gamma(t) \, \mathrm{d}\gamma  \leq \int_0^\infty \frac{\gamma^{\delta-1}l(t)}{1+\gamma(1\ast l)(t)} \, \mathrm{d}\gamma = \Gamma(1-\delta) l(t)[(1\ast l)(t)]^{-\delta},\quad t>0,
			\end{align*}
			
			Finally, simple calculation can derives that
			\begin{align*}
				&\widehat{c}_\delta(\lambda) =- \int_0^\infty \, \theta^{-\delta}\widehat{\phi}(\lambda,\mathrm{d}\theta)=(\lambda \widehat{l}(\lambda))^{-1}\int_0^\infty \theta^{-\delta} e^{-\theta \widehat{l}(\lambda)^{-1}}\mathrm{d}\theta = \Gamma(1-\delta)\lambda^{-1}\widehat{l}(\lambda)^{-\delta},\\
				&\widehat{d}_\delta(\lambda) = \int_0^\infty \, \theta^{-\delta}\widehat{\eta}(\lambda,\mathrm{d}\theta)=\int_0^\infty \theta^{-\delta} e^{-\theta \widehat{l}(\lambda)^{-1}}\mathrm{d}\theta = \Gamma(1-\delta)\widehat{l}(\lambda)^{1-\delta}.
			\end{align*}
			We complete the proof.
		\end{proof}
	\end{lemma}
	
	\begin{remark}
		For each $\delta \in (0,1)$, the function $c_\delta(t)$ defined in Lemma \ref{A.2} additionally satisfies the following estimates:
		\begin{align*}
			c_\delta(t) \geq \Gamma(1-\delta) k_1(t)^{\delta }, \quad t>0.
		\end{align*}
		In fact, the left-hand side inequality in Proposition \ref{prop 2.1}(ii) alone suffices to yield
		\begin{align*}
			\frac{1}{\Gamma(\delta)}\int_0^\infty \gamma^{\delta-1} s_\gamma(t) \, \mathrm{d}\gamma  \geq \frac{1}{\Gamma(\delta)}\int_0^\infty \frac{\gamma^{\delta-1}}{1+\gamma [k_1(t)]^{-1}} \, \mathrm{d}\gamma = \Gamma(1-\delta) k_1(t)^{\delta},\quad t>0.
		\end{align*}
	\end{remark}
	
	\section{New completely positive kernels}\label{app b}
	The following lemma provides a novel approach to inducing completely positive kernels and serves as our fundamental tool for handling fully nonlocal telegraph equations.
	\begin{lemma}\label{lemma B.1}
		Assume that $l\in L^1_{loc}(\mathbb R_+)$ is completely positive.
		Let $\delta_1, \delta_2 \in (0,1)$ with $\delta_1+\delta_2 \geq 1$, $\gamma \in \mathbb R_+$, and $r_\gamma(t)$ be defined by \eqref{2.4}. Then there exists a completely positive function $h_{\delta_1,\delta_2, \gamma}:(0,\infty) \to (0,\infty)$, whose Laplace transform is given by
		\begin{align}\label{B.1}
			\widehat{h}_{\delta_1,\delta_2,\gamma}(\lambda):=\Gamma(1-\delta_1)\Gamma(1-\delta_2) \widehat{l}(\lambda)^{1-\delta_1}\widehat{r}_\gamma(\lambda)^{1-\delta_2}.
		\end{align}
		Moreover, if $b$ is also nonincreasing, then there exists $C_{\delta_1,\delta_2}>0$ depends solely on $\delta_1, \delta_2$, such that
		\begin{align*}
			h_{\delta_1,\delta_2,\gamma}(t) \leq C_{\delta_1,\delta_2} l\left(\frac{t}{2}\right)\left[(1 \ast l)\left(\frac{t}{2}\right)\right]^{1-(\delta_1+\delta_2)}, \quad t>0.
		\end{align*}
		\begin{proof}
			We deduce from Proposition \ref{prop 2.1}(i) that
			\begin{align*}
				k_0 r_\gamma(t)+((k_1 +\gamma 1)\ast r_\gamma)(t) = 1, \quad t>0, \ \gamma \in \mathbb R_+.
			\end{align*}
			Since $k_0\geq 0$ and $k_1+\gamma 1$ is nonnegative and nonincreasing, $r_\gamma$ is completely positive. By Lemma \ref{lemma A.1} and Lemma \ref{lemma A.2}, for each $\gamma \in \mathbb R_+$, there exists a nonnegative function $d_{\delta_2,\gamma} \in L^1_{loc}(\mathbb R_+)$ admits the representation
			\begin{align*}
				d_{\delta_2,\gamma}(t):= \int_0^\infty \theta^{-\delta_2} \eta_\gamma(t, \mathrm{d} \theta),  \quad t>0,
			\end{align*}
			whose Laplace transform is $\widehat{d}_{\delta_2,\gamma}(\lambda) = \Gamma(1-\delta_2)\widehat{r}_\gamma(\lambda)^{1-\delta_2}$, where $\eta_\gamma(t,\mathrm{d}\theta)$ is a positive finite measure such that
			$\widehat{\eta}_\gamma(\lambda, \mathrm{d}\theta) = e^{-\theta \widehat{r}_\gamma(\lambda)^{-1}}\mathrm{d}\theta$.
			
			Since $d_{\delta_1}, d_{\delta_2,\gamma} \in L_{loc}^1(\mathbb R_+)$, \cite[Theorem 2.2(i)]{G. Gripenberg} shows that $h_{\delta_1,\delta_2,\gamma}:=(d_{\delta_1} \ast d_{\delta_2,\gamma})(t)$ is well-defined for $t>0$ and $h_{\delta_1,\delta_2,\gamma} \in L_{loc}^1(\mathbb R_+)$. Let us to prove that $h_{\delta_1,\delta_2,\gamma}$ is completely positive.
			Applying the Laplace transform to
			$h_{\delta_1,\delta_2,\gamma}$ yields
			\begin{align*}
				\widehat{h}_{\delta_1,\delta_2,\gamma}(\lambda)=\widehat{d}_{\delta_1}(\lambda) \widehat{d}_{\delta_2,\gamma}(\lambda)=\Gamma(1-\delta_1)\Gamma(1-\delta_2)\widehat{l}(\lambda)^{1-\delta_1} \widehat{r}_\gamma(\lambda) ^{1-\delta_2}.
			\end{align*}
			Note that $\widehat{r}_\gamma(\lambda)=\widehat{l}(\lambda)(1+\gamma \widehat{l}(\lambda))^{-1}$. Substituting this into the above, we obtain
			\begin{align}\label{B.2}
				\widehat{h}_{\delta_1,\delta_2,\gamma}(\lambda)^{-1} =\Gamma(1-\delta_1)\Gamma(1-\delta_2) \left[\widehat{l}(\lambda)^{-1}\right]^{1-\delta_1} \left[\widehat{l}(\lambda)^{-1}+\gamma \right]^{1-\delta_2}.
			\end{align}
			Due to the assumption of $\delta_1, \delta_2$ and \cite[Proposition 7.13]{R. L. Schilling} we obtain that $\lambda \mapsto \lambda^{1-\delta_1}(\lambda+\gamma)^{1-\delta_2}$
			is a Bernstein function. Meanwhile, since the kernel $l$ is completely positive, \cite[Proposition 4.5]{J. Pruss} yields that $[\widehat{l}(\lambda)]^{-1}$ is a Bernstein function. Note that $\mathcal{BF}\circ\mathcal{BF} \in \mathcal{BF}$(\cite[Proposition 4.2]{J. Pruss}). Then \eqref{B.2} shows that $\widehat{h}_{\delta_1,\delta_2,\gamma}(\lambda)^{-1} \in (\mathcal BF)$, which
			implies $h_{\delta_1,\delta_2,\gamma}$ is completely positive by \cite[proposition 4.5]{J. Pruss}.
			
			Let us consider the second assertion. If $l$ is nonincreasing, then from Proposition \ref{prop 2.1}(iii) we get
			\begin{align*}
				r_\gamma(t)[(1\ast r_\gamma)(t)]^{-\delta_2}&\leq
				\frac{l(t)}{1+\gamma(1\ast l)(t)}\left[\frac{(1\ast l)(t)}{1+\gamma (1\ast l)(t)}\right]^{-\delta_2}\\
				&= l(t)[(1\ast l)(t)]^{-\delta_2}[1+\gamma(1\ast l)(t)]^{\delta_2-1}\\
				&\leq l(t)[(1\ast l)(t)]^{-\delta_2}, \quad t>0.
			\end{align*}
			Combining the preceding inequality with Lemma \ref{lemma A.2}, we deduce that
			\begin{align*}
				h_{\delta_1,\delta_2,\gamma}(t)&= \int_0^t d_{\delta_1}(t-\zeta) d_{\delta_2,\gamma}(\zeta) \, \mathrm{d}\zeta\\
				&\lesssim \int_0^t l(t-\zeta)[(1\ast l)(t-\zeta)]^{-\delta_1} r_\gamma(\zeta)[(1\ast r_\gamma)(\zeta)]^{-\delta_2} \, \mathrm{d}\zeta\\
				&\lesssim \int_0^t l(t-\zeta)[(1\ast l)(t-\zeta)]^{-\delta_1} l(\zeta)[(1\ast l)(\zeta)]^{-\delta_2}\, \mathrm{d}\zeta\\
				&\lesssim l\left(\frac{t}{2}\right)\left[(1\ast l)\left(\frac{t}{2}\right)\right]^{-\delta_1} \int_0^{\frac{t}{2}}  l(\zeta) [(1\ast l)(\zeta)]^{-\delta_2} \, \mathrm{d}\zeta\\
				&\quad + l\left(\frac{t}{2}\right)\left[(1\ast l)\left(\frac{t}{2}\right)\right]^{-\delta_2} \int_{\frac{t}{2}}^t  l(t-\zeta)[(1\ast l)(t-\zeta)]^{-\delta_1}\, \mathrm{d}\zeta\\
				&\lesssim  l\left(\frac{t}{2}\right)\left[(1 \ast l)\left(\frac{t}{2}\right)\right]^{1-(\delta_1+\delta_2)}.
			\end{align*}
			Take $C_{\delta_1,\delta_2}=\frac{(2-\delta_1-\delta_2)\Gamma(1-\delta_1)\Gamma(1-\delta_2)}{(1-\delta_1)(1-\delta_2)}$. The proof is completed.
		\end{proof}
	\end{lemma}
	We term $h_{\gamma,\delta_1,\delta_2}$ the subordinate completely positive kernel associated with the completely positive kernel $l$. In view of \cite[Theorem 2.2]{P. Clement 81}, there exist $m_0\geq0$ and a nonnegative nonincreasing function $m_1 \in L^1_{loc}(\mathbb R_+)$ such that
	\begin{align}\label{B.3}
		m_0h_{\delta_1,\delta_2,\gamma} + (m_1\ast h_{\delta_1,\delta_2,\gamma})(t) =1 , \quad t>0.
	\end{align}
	We maintain the decomposition $m_1(t)=m_\infty +m_2(t)$, where $m_2\in L^1_{loc}(\mathbb R_+)$ is nonnegative and nonincreasing, $m_\infty = \lim\limits_{t \to \infty}m_1(t)$. The next lemma gives the expression of $m_0$ and the asymptotic behavior of $m_1(t)$ at $t=\infty$.
	\begin{lemma}\label{lemma B.2}
		Assume that $l\in L^1_{loc}(\mathbb R_+)$
		is completely positive. Let $\delta_1, \delta_2 \in (0,1)$ with $\delta_1+\delta_2 \geq 1$, then the following hold:
		\begin{itemize}
			\item[(i)] $m_\infty =\frac{ k_\infty^{1-\delta_1}(k_\infty+\gamma)^{1-\delta_2}}{\Gamma(1-\delta_1)\Gamma(1-\delta_2) }$.
			\item[(ii)] $m_0 =0$ if $\delta_1+\delta_2>1$, and $m_0=\frac {k_0}{ \Gamma(1-\delta_1)\Gamma(1-\delta_2)}$ if $\delta_1+\delta_2 =1$.
		\end{itemize}
		\begin{proof}
			It follows from \eqref{B.1} and \eqref{B.3} that
			\begin{align*}
				\frac{\widehat{l}(\lambda)^{\delta_1-1} \left[\widehat{l}(\lambda)^{-1}+\gamma \right]^{1-\delta_2} }{\Gamma(1-\delta_1)\Gamma(1-\delta_2)}= \lambda (m_0 + \frac{m_\infty}{\lambda} + \widehat{m_2}(\lambda)),
			\end{align*}
			Substituting $\widehat{l}(\lambda)^{-1} = \lambda (k_0+ \frac{k_\infty}{\lambda}+\widehat{k_2}(\lambda))$ into the above yields
			\begin{align}\label{B.4}
				\frac{\left[k_0 \lambda+ k_\infty+\lambda\widehat{k_2}(\lambda)\right]^{1-\delta_1} \left[k_0\lambda+k_\infty +\lambda\widehat{k_2}(\lambda)+\gamma \right]^{1-\delta_2} }{\Gamma(1-\delta_1)\Gamma(1-\delta_2)} = \lambda m_0 + m_\infty + \lambda\widehat{m_2}(\lambda),
			\end{align}
			Taking the limit as $\lambda \to 0$ on both sides of the above gives
			\begin{align*}
				m_\infty =k_\infty^{1-\delta_1}(k_\infty+\gamma)^{1-\delta_2}.
			\end{align*}
			On the other hand, note that \eqref{B.4} is equivalent to
			\begin{align*}
				\frac{\lambda^{1-(\delta_1+\delta_2)}\left[k_0+\frac{k_\infty}{\lambda}+\widehat{k_2}(\lambda)\right]^{1-\delta_1} \left[k_0+ \frac{k_\infty}{\lambda}+\widehat{k_2}(\lambda)+\frac{\gamma}{\lambda} \right]^{1-\delta_2}}{\Gamma(1-\delta_1)\Gamma(1-\delta_2)}= m_0 + \frac{m_\infty}{\lambda} + \widehat{m_2}(\lambda).
			\end{align*}
			Taking the limit as $\lambda \to \infty$ on both sides of the above. This completes the proof.
		\end{proof}
	\end{lemma}
	
	\section{Estimates for the fractional Laplacian}\label{app c}
	Early research on the decay estimates for the operator $\mathcal F^{-1}(|\xi|^{-s}e^{it|\xi|^{\alpha}}\mathcal F )$ with $s>0$ can be traced back to \cite{H. Pecher76} by Pecher. He considered the case where $\alpha$ is an integer. In fact,
	employing the approach of Pecher, we can readily generalize to derive the case where $\alpha \in (\frac{1}{2},1)$. For other related research, see \cite{H. Pecher96, M. Sugimoto, A. Miyachi, M. DAbbicco,M. DAbbicco14}. 
	
	\begin{lemma}\label{lemma C.1}
			Let $s \in (0,3)$. Assume that $\beta, p$ satisfy
		\begin{align}\label{C.1}
			\begin{cases}
				p \in (1, \frac{6}{3+s}] \cap [\frac{4}{2+s}, \frac{6}{3+s}] \text{ ~if~ } \beta=2,\\
				p \in (1, \frac{6}{3+s}] \cap [ 2\frac{12-3\beta}{12-3\beta+ +4s}, \frac{6}{3+s}] \text{ ~if~ }\beta \in (1,2).
			\end{cases}
		\end{align}
		It holds that
		\begin{align}\label{C.2}
			\left\|\mathcal F^{-1} \left(|\xi|^{-s} e^{\pm i t|\xi|^{\frac{\beta}{2}}} \tilde{u}(\xi) \right) \right\|_{L^{p'}} \lesssim t^{\frac{2s}{\beta}-\frac{6}{\beta}(\frac{1}{p} - \frac{1}{p'} )} \|u\|_{L^{p}} \text{ ~for~ } t > 0.
		\end{align}
		\begin{proof}
			The case $ \beta = 2 $ has been established in \cite[Theorem 2.2]{H. Pecher76}(cf. also \cite[Theorem 1.1]{H. Pecher96}). It therefore suffices to consider the scenario where $ \beta \in (1,2) $. Define $Q(\xi):=|\xi|^{\frac{\beta}{2}}$ with $\xi \in \mathbb R^3$. A straightforward computation yields the Hessian matrix
			\begin{align*}
				H_Q =\frac{\beta}{2}r^{\frac{\beta}{2} - 4}(r^2 I+(\frac{\beta}{2}  -2) \xi \xi^T )  \text{ ~for~ } \xi \neq 0,
			\end{align*}
			where $I$ is the identity matrix,  $\xi^T$ denotes the transpose of
			$\xi$ and $r = \sqrt{\xi_1^2 +\xi_2^2 +\xi_3^2}$. Then, via standard linear algebra steps, it can be derived that the eigenvalues of $H_Q$ are
			\begin{align*}
				\frac{\beta}{2}(\frac{\beta}{2} -1)r^{\frac{\beta}{2}-2}, \ \frac{\beta}{2}r^{\frac{\beta}{2}-2 }(\text{double}) \text{ ~for~ } \xi \neq 0.
			\end{align*}
			Note that $H_Q$ is a real symmetric matrix. then
			\begin{align*}
				\text{rank}(H_Q) =
					3 \text{ ~for~ } \beta \in (1,2) , \ \xi \neq 0.
			\end{align*}
			
			Denote $Q_1(\xi) = |\xi|^{s}$ with $\xi \in \mathbb R^3$. We now proceed to prove
			\begin{align}\label{C.3}
				\left\|\mathcal F^{-1}\left(|\xi|^{-s } e^{\pm i |\xi|^{\frac{\beta}{2} }} \varrho_j(\xi) \tilde{u}(\xi) \right)\right \|_{L^{\infty}}
				\lesssim 2^{(3-s- \frac{ 3}{4}\beta) j}\| u\|_{L^1},  \quad j \in \mathbb N_+.
			\end{align}
			Due to $\varrho/Q_1 \in C_c^\infty$, invoking Lemma \ref{lemma 2.1} yields
			\begin{align}\label{C.4}
				\begin{aligned}
					\left\|\mathcal F^{-1}\left(|\xi|^{-s }e^{\pm i 2^{ \frac{\beta j}{2} }Q(\xi)}  \varrho(\xi) \right)\right\|_{L^\infty}
					\leq L(\varrho / Q_1, Q) 2^{-\frac{3\beta j}{4}} \text{ ~for~ } j \in \mathbb N_+,
				\end{aligned}
			\end{align}
			where $L(\varrho/Q_1, Q)$  is a constant that depends solely on $\varrho/Q_1, Q$. Note that
			\begin{align*}
				\left\|\mathcal F^{-1}\left(|\xi|^{-s } e^{\pm i |\xi|^{\frac{\beta}{2} }} \varrho_j(\xi) \right)\right\|_{L^\infty} = 2^{(3-s)j} \left\|\mathcal F^{-1}\left(|\xi|^{-s } e^{\pm i 2^{ \frac{\beta j}{2} }|\xi|^{\frac{\beta}{2} }} \varrho(\xi) \right)\right\|_{L^\infty}, \quad j \in \mathbb N_+.
			\end{align*}
			This along with \eqref{C.4} shows that
			\begin{align*}
				\left\|\mathcal F^{-1}\left(|\xi|^{-s } e^{ \pm i |\xi|^{\frac{\beta}{2} } } \varrho_j(\xi) \right)\right\|_{L^\infty}\lesssim 2^{(3-s- \frac{3}{4}\beta) j}, \quad j \in \mathbb N_+.
			\end{align*}
			It follows from the above inequality and the Young' inequality that \eqref{C.3} holds. Next, let us  establish the following:
			\begin{align}\label{C.5}
				\left\|\mathcal F^{-1}\left(|\xi|^{-s} e^{\pm i  |\xi|^{\frac{\beta}{2} }} \varrho_j(\xi) \tilde{u}(\xi) \right)\right \|_{L^{2}}
				\lesssim 2^{-{s j}}\| u\|_{L^2},  \quad j \in \mathbb N_+.
			\end{align}
			One can readily observe that
			\begin{align*}
				\left \||\xi|^{-s } e^{ \pm i|\xi|^{\frac{\beta}{2} } } \varrho_j(\xi)\right \|_{L^\infty} \lesssim 2^{-{s j}} \text{ ~for~ } j \in \mathbb N^+.
			\end{align*}
			This along with the Plancherel's identity shows that \eqref{C.5} holds.
			Applying the Riesz-Thorin interpolation theorem to \eqref{C.3} and \eqref{C.5}, we conclude that
			\begin{align}\label{C.6}
				\left \|\mathcal F^{-1}\left(|\xi|^{-s } e^{\pm i |\xi|^{\frac{\beta}{2} }} \varrho_j(\xi) \tilde{u}(\xi) \right)\right \|_{L^{p'}} \lesssim 2^{((3-\frac{3\beta}{4})(\frac{1}{p}- \frac{1}{p'} ) -s )j}\|u\|_{L^p} \leq\|u\|_{L^p}, \quad j \in \mathbb N_+.
			\end{align}
			Here, we use the fact $(3-\frac{3\beta}{4})(\frac{1}{p}- \frac{1}{p'} ) -s\leq 0$,
			provided the assumption \eqref{C.1}.
			On the other hand, note that $1<p \leq \frac{6}{3+s}$ implies $\frac{1}{p}-\frac{1}{p'}\in [\frac{s}{3},1)$. Let $\psi$ satisfy $\frac{1}{\psi}=\frac{1}{p}-\frac{1}{p'}$, then $\psi \in (1, \frac{3}{s} ]$. Since $\varrho_0 \in C_c^\infty(\mathbb R^3)$, it follows directly that
			\begin{align*}
				|\{\xi \in \mathbb R^3: ||\xi|^{-s} e^{\pm i |\xi|^{\frac{\beta}{2} } }\varrho_0(\xi) | \geq \alpha \}| \lesssim \min\{1, \alpha^{-\frac{3}{s}}\} \lesssim \alpha^{-\psi}.
			\end{align*}
			Invoking Lemma \ref{lemma 2.2} yields
			\begin{align}\label{C.7}
				\left \|\mathcal F^{-1}\left(|\xi|^{-s } e^{\pm i |\xi|^{\frac{\beta}{2} }} \varrho_0(\xi) \tilde{u}(\xi) \right)\right \|_{L^{p'}} \lesssim \|u\|_{L^p}.
			\end{align}
			
			The next step consists in deriving from \eqref{C.6} and \eqref{C.7} that:
			\begin{align}\label{C.8}
				\left\| \mathcal F^{-1} \left(|\xi|^{-s } e^{ \pm i |\xi|^{\frac{\beta}{2} }}  \tilde{u}(\xi) \right) \right\|_{B^0_{p',2}} \lesssim \left\|u \right\|_{B^0_{p,2}}.
			\end{align}
			Note that $\mathcal{F}^{-1}(\varrho_j) \ast u = \sum\limits_{j_1=j-1}^{j+1 }\mathcal{F}^{-1} (\varrho_j) \ast \mathcal F^{-1}(\varrho_{j_1}) \ast u$
			for $j \in \mathbb N$ with $\varrho_{-1} = 0$. Substituting
			$\varrho_j \ast u$ for $u$ in
			\eqref{C.6} and \eqref{C.7} yields
			\begin{align}\label{C.9}
				\begin{aligned}
					&\quad \left\| \mathcal F^{-1} \left(|\xi|^{-s } e^{ \pm i |\xi|^{\frac{\beta}{2} }}  \varrho_j(\xi)\tilde{u}(\xi) \right) \right\|_{L_{p'}}\\
					&\leq \sum_{j_1 =j-1 }^{j+1 } \left\| \mathcal F^{-1} \left(|\xi|^{-s } e^{ \pm i |\xi|^{\frac{\beta}{2} }}  \varrho_j(\xi)\varrho_{j_1}(\xi)\tilde{u}(\xi) \right) \right\|_{L_{p'}} \\
					&\lesssim \sum_{j_1 =j-1 }^{j+1 } \left\| \mathcal F^{-1} \left(|\xi|^{-s } e^{ \pm i |\xi|^{\frac{\beta}{2} }}  \varrho_{j_1}(\xi)\tilde{u}(\xi) \right) \right\|_{L_{p}} \text{ ~for~ } j \in \mathbb N.
				\end{aligned}
			\end{align}
			Then raising both sides of \eqref{C.9} to the second power and summing over all indices $j$ yields \eqref{C.8} as a consequence. Note that the embeddings
			$B^0_{p',2} \hookrightarrow L_{p'}$ and $L_p \hookrightarrow B^0_{p,2}$. We immediately have
			\begin{align*}
				\left\|\mathcal F^{-1} \left(|\xi|^{-s} e^{\pm i |\xi|^{\frac{\beta}{2}}} \tilde{u}(\xi) \right) \right\|_{L^{p'}} \lesssim \|u\|_{L^{p}}.
			\end{align*}
			Finally, by homogeneity of $Q(\xi)$ we obtain
			\begin{align*}
				\begin{aligned}
					&\quad\left\|\mathcal F^{-1} \left(|\xi|^{-s} e^{ \pm i t|\xi|^{\frac{\beta}{2}}} \tilde{u}(\xi) \right)(x)\right\|_{L^{p'}}\\
					&=t^{\frac{2s}{\beta} } \left\|\mathcal F^{-1} \left(|\xi|^{-s} e^{\pm i |\xi|^{\frac{\beta}{2}}} \tilde{u}_{t^{\frac{2}{\beta}}}(\xi) \right)\left(t^{-\frac{2}{\beta} } x \right) \right\|_{L^{p'}}\\
					&\lesssim t^{\frac{2s}{\beta}  +\frac{6}{\beta p'}}  \left\|\mathcal F^{-1} \left(|\xi|^{-s} e^{\pm i |\xi|^{\frac{\beta}{2}}} \tilde{u}_{t^{\frac{2}{\beta}}}(\xi) \right)(x) \right\|_{L^{p'}}\\
					&\lesssim t^{\frac{2s}{\beta} +\frac{6}{\beta p'}}  \left\| u_{t^{\frac{2}{\beta}}}(x)\right\|_{L^{p}}\\
					&\lesssim t^{\frac{2s}{\beta}-\frac{6}{\beta}(\frac{1}{p} - \frac{1}{p'} )}\|u\|_{L^p},
				\end{aligned}
			\end{align*}
			where $u_{t^{\frac{2}{\beta}}}(x)= u(t^{\frac{2}{\beta}}x)$. We complete the proof.
		\end{proof}
	\end{lemma}
	We can further derive the following practical estimates.
	\begin{proposition}\label{prop C.1}
		Let $s \in (0,3)$. Assume that $\beta, p$ satisfy \eqref{C.1}. It holds that
		\begin{align*}
			\left\|\mathcal F^{-1} \left(\cos( t|\xi|^{\frac{\beta}{2}}) \tilde{u}(\xi) \right) \right\|_{L^{p'}} \lesssim |t|^{\frac{2s}{\beta} - \frac{6}{\beta}(\frac{1}{p} - \frac{1}{p'})}\|u\|_{H^{s,p}} \text{ ~for~ } t \neq 0.
		\end{align*}
		\begin{proof}
			Without loss of generality, we may restrict our consideration to the case $t>0$. Lemma \ref{C.1} then yields
			\begin{align*}
				&\quad \left\|\mathcal F^{-1}\left(\cos (t |\xi|^{\frac{\beta}{2} })\tilde u(\xi) \right)\right\|_{L^{p'}} \\
				&\leq \left\|\mathcal F^{-1}\left(|\xi|^{-s }e^{ it |\xi|^{\frac{\beta}{2} }}  \tilde |\xi|^s \tilde u(\xi) \right)\right\|_{L^{p'}}+\left\|\mathcal F^{-1}\left(|\xi|^{-s } e^{-it |\xi|^{\frac{\beta}{2} }\tilde u(\xi)}  \tilde u(\xi) \right)\right\|_{L^{p'}}\\
				&\lesssim t^{\frac{2s}{\beta}-\frac{6}{\beta}(\frac{1}{p} - \frac{1}{p'} )}\|\mathcal F^{-1}(|\xi|^{s}\tilde u(\xi))\|_{L^p}\\
				&\lesssim t^{\frac{2s}{\beta}-\frac{6}{\beta}(\frac{1}{p} - \frac{1}{p'} )}\| u)\|_{H^{s,p}}.
			\end{align*}
			This complete the proof.
		\end{proof}
	\end{proposition}
	
	\begin{proposition}\label{prop C.2}
		Assume that $\beta, p$ satisfy
		\begin{align}\label{C.10}
			\begin{cases}
				p \in [\frac{4}{3}, 2] \text{ ~if~ } \beta =2,\\
				p \in [2\frac{12-3\beta}{12-\beta}, 2] \text{ ~if~ }\beta \in (1,2).
			\end{cases}
		\end{align}
		It holds that
		\begin{align}\label{C.11}
			\left\|\mathcal F^{-1} \left(|\xi|^{-\frac{\beta}{2}} \sin( t|\xi|^{\frac{\beta}{2}}) \tilde{u}(\xi) \right) \right\|_{L^{p'}} \lesssim |t|^{1 - \frac{6}{\beta}(\frac{1}{p} - \frac{1}{p'})}\|u\|_{L^{p}} \text{ ~for~ } t \neq 0.
		\end{align}
		\begin{proof}
			Without loss of generality, we may restrict our consideration to the case $t>0$. It is straightforward to find that
			\begin{align*}
				\left\|\mathcal F^{-1}\left(|\xi|^{-\frac{\beta}{2}} \sin(  |\xi|^{\frac{\beta}{2}})  \varrho_0(\xi) \right)\right\|_{L^\infty}< \infty, \  \left\| |\xi|^{-\frac{\beta}{2}} \sin(  |\xi|^{\frac{\beta}{2}})  \varrho_0(\xi) \right\|_{L^\infty}< \infty.
			\end{align*}
			These along with the Young's inequality implies
			\begin{align}\label{C.12}
				\begin{aligned}
					&\left\|\mathcal F^{-1}\left(|\xi|^{-\frac{\beta}{2} } \sin(  |\xi|^{\frac{\beta}{2} }) \varrho_0(\xi) \tilde{u}(\xi) \right)\right \|_{L^{\infty}}
					\lesssim \| u\|_{L^1}, \\
					&\left\|\mathcal F^{-1}\left(|\xi|^{-\frac{\beta}{2} } \sin(  |\xi|^{\frac{\beta}{2} }) \varrho_0(\xi) \tilde{u}(\xi) \right)\right \|_{L^{2}}
					\lesssim \| u\|_{L^2},
				\end{aligned}
			\end{align}
			respectively. Applying the Riesz-Thorin interpolation theorem to \eqref{C.12} again, one can derive that
			\begin{align*}
				\left \|\mathcal F^{-1}\left(|\xi|^{-\frac{\beta}{2} } \sin( |\xi|^{\frac{\beta}{2} }) \varrho_0(\xi) \tilde{u}(\xi) \right)\right \|_{L^{p'}} \lesssim \|u\|_{L^p}.
			\end{align*}
			Obviously, we have
			\begin{align*}
				&\quad \left\|\mathcal F^{-1}\left(|\xi|^{-\frac{\beta}{2} } \sin(|\xi|^{\frac{\beta}{2} })\varrho_j(\xi) \right)\right\|_{L^{p'}}\\
				&\leq \left\|\mathcal F^{-1}\left(|\xi|^{-\frac{\beta}{2} }e^{i |\xi|^{\frac{\beta}{2} }}  \varrho_j(\xi) \right)\right\|_{L^{p'}}+\left\|\mathcal F^{-1}\left(|\xi|^{-\frac{\beta}{2} } e^{-i |\xi|^{\frac{\beta}{2} }}  \varrho_j(\xi) \right)\right\|_{L^{p'}} \text{ ~for~ } j \in \mathbb N_+
			\end{align*}
			which implies
			\begin{align*}
				\left\|\mathcal F^{-1}\left(|\xi|^{-\frac{\beta}{2} } \sin(|\xi|^{\frac{\beta}{2} })\varrho_j(\xi) \right)\right\|_{L^{p'}} \lesssim 2^{((3-\frac{3\beta}{4})(\frac{1}{p}- \frac{1}{p'} ) -\frac{\beta}{2} )j}\|u\|_{L^p} \leq\|u\|_{L^p}, \quad j \in \mathbb N_+
			\end{align*}
			by \eqref{C.6}.
			Here, we use the fact that
			$(3-\frac{3\beta}{4})(\frac{1}{p}- \frac{1}{p'} ) -\frac{\beta}{2}\leq 0$,
			provided the assumption \eqref{C.10}. The following proof follows the same argument as in the proof subsequent to \eqref{3.7} in Lemma \ref{lemma C.1}, and is therefore omitted here.
		\end{proof}
	\end{proposition}
	
	\begin{remark}
		For $\beta=2$, Proposition \ref{prop C.2} reduces to a special case of \cite[Theorem]{P. Brenner}.
	\end{remark}
	\end{appendices}
	
	\section*{Statements and Declarations}
	\noindent{\bf Competing interests}\\
	The authors declare that there is no conflict of interests regarding the publication of this paper.\\
	\noindent{\bf Funding}\\
	This work was supported by National Natural Science Foundation of China (12471172) and the Project of Hunan Provincial Education Department of China (23B0125,24B0168).\\
	\noindent{\bf Availability of data and materials}\\
	Not applicable.
	
	\renewcommand{\thesection}{\arabic{section}}
	
\end{document}